\documentclass[10pt]{article}
\usepackage{graphicx}
\usepackage{amssymb}
\usepackage{epstopdf}
\usepackage{amsmath}
\usepackage{amsfonts}
\newcommand{\braket}[2]{\langle #1,#2 \rangle}
\newcommand{\la}{\lambda}

\def\phi{{\varphi}}

\DeclareSymbolFont{AMSb}{U}{msb}{m}{n}
\DeclareMathSymbol{\N}{\mathbin}{AMSb}{"4E}
\DeclareMathSymbol{\Z}{\mathbin}{AMSb}{"5A}
\DeclareMathSymbol{\R}{\mathbin}{AMSb}{"52}
\DeclareMathSymbol{\Q}{\mathbin}{AMSb}{"51}
\DeclareMathSymbol{\I}{\mathbin}{AMSb}{"49}
\DeclareMathSymbol{\C}{\mathbin}{AMSb}{"43}

\begin{document}

\addtolength{\textheight}{0 cm} \addtolength{\hoffset}{0 cm}

\def\e{\varepsilon}

\def\bR{\mathbb{R}}

\def\cA{ \mathcal{A} }
\def\cF{ \mathcal{F} }
\def\cK{ \mathcal{K} }
\def\cL{ \mathcal{L} }
\def\cT{\mathcal{T}}

\def\dL{ \overline{\partial } L }
\def\dM{ \overline{\partial } M }
\def\dF{ \overline{\partial } F }

\def\va{ \vec{a}}

\def\S{ \int_{\Omega} }
\def\dx{\,\mathrm{d}x}
\def\dy{\,\mathrm{d}y}
\def\dt{\,\mathrm{d}t}

\def\Lp{L^p(\Omega)}
\def\Ldos{L^2(\Omega)}
\def\Huno{H_0^1(\Omega)}

\newtheorem{thm}{Theorem}[section]
\newtheorem{pro}{Proposition}[section]
\newtheorem{lem}{Lemma}[section]
\newtheorem{cor}{Corollary}[section]
\newtheorem{defn}{Definition}[section]
\setcounter{secnumdepth}{5}
\newtheorem{remark}{Remark}[section]
\newtheorem{application}{Application}

\title{Homogenization of maximal monotone vector fields via selfdual variational calculus}

\author{ \ Nassif Ghoussoub \thanks{Partially supported by a grant
from the Natural Sciences and Engineering Research Council of Canada.}\\
\small Department of Mathematics,
\small University of British Columbia, \\
\small Vancouver BC Canada V6T 1Z2 \\
\small {\tt nassif@math.ubc.ca}\\
\and \ Abbas Moameni\thanks{Research supported by a Coleman  fellowship at Queen's University.} \\
\small Department of Mathematics and Statistics,
\small Queen's University, \\
\small Kingston ON  Canada K7L 3N6 \\
 \small {\tt momeni@mast.queensu.ca}\\
\and \ Ram\'on Z\'arate S\'aiz \thanks{Partially supported by a graduate fellowship from CONACYT (Mexico) and is part of a PhD thesis under the supervision of N. Ghoussoub}\\
\small Department of Mathematics,
\small University of British Columbia, \\
\small Vancouver BC Canada V6T 1Z2 \\
\small {\tt ramon@math.ubc.ca}
\\ \date{} }

\maketitle

\section*{Abstract}
We use the theory of selfdual Lagrangians to give a variational approach to the homogenization of equations in divergence form, that are driven by a periodic family of maximal monotone vector fields. The approach has the advantage of using $\Gamma$-convergence methods for corresponding  functionals just as in the classical case of convex potentials, as opposed to the graph convergence methods used in the absence of potentials. 
A new variational formulation for the homogenized equation is also given. 

\tableofcontents
%\newpage
  \section{Introduction}
  We consider the homogenization of the problem 
  \begin{eqnarray}\label{eleven}
\left\{
  \begin{array}{lll}
 \hfill   \tau_n(x) &\in\quad  {\ss}(\frac{x}{\e_n}, \nabla u_n(x)) &x \in \Omega, \\
   -{\rm div} (\tau_n (x))&=\quad u^*_n(x) & x \in \Omega,\\
\hfill u_n(x)&=\quad 0 & x\in \partial \Omega,\\
  \end{array}
\right.
\end{eqnarray}
where $\Omega$ is a bounded  domain of $\bR^N$  and  ${\ss}: \Omega \times \bR^N \to \bR^N$  is a measurable map on $\Omega \times \bR^N$ such that $\ss (x, \cdot )$ is maximal monotone on $\bR^N$ for almost all $x\in \Omega$, and such that  ${\ss}(., \xi)$ is  $Q$-periodic for an  open non-degenerate parallelogram $Q$ in  $\bR^n$.  This problem has been investigated  in recent years by  many authors. We refer the interested reader to \cite{Att-book, Br-Ch, Br-Ch-D, Ch-dalM-D, D-M-S, Fr-Mu-Ta, FM1, PD} for related results.

The particular case where the maximal monotone operator is a subdifferential of the form 
\begin{equation}\label{twelve}
{\ss}(x, \xi)=\partial_{\xi}\psi(x, \xi),
\end{equation}
with $\psi: \Omega \times \bR^N \to \bR $ being a convex  function in the second variable is particularly appealing and completely  understood. Indeed, under appropriate boundedness and coercivity conditions on $\psi$, say 
\[C_0(|\xi|^p-1)\leq \psi(x, \xi) \le C_1(|\xi|^p+1) \qquad \text{ for all } (x, \xi) \in \Omega \times \bR^N, \]
where $1<p<\infty$   and $C_0, C_1$ are positive constants, 
one can then use a variational approach to identify for a given $u^* \in W^{-1,p}(\Omega)$, the solution $(u, \tau)$ of (\ref{eleven}) as the respective minima of the problems 
\begin{equation}\label{thirteen}
\inf\left\{ \int_\Omega \psi (x, \nabla u (x))\, dx - \int_\Omega u^*(x) u(x)\, dx; \,\, u\in  W_0^{1,p}(\Omega)\right\},
\end{equation}
 and  
\begin{equation}
\inf\left\{  \int_\Omega \psi^* (x, \tau (x))\, dx;\,\, {\rm div}(\tau)=u^*\right\},    
\end{equation}
where $\psi^*$ is the Fenchel-Legendre dual (in the second variable) of $\psi$. In this case, the classical concept of $\Gamma$-convergence --introduced by DeGiorgi-- can be used to show that if 
 $u^*_n \to u^* $ strongly in $W^{-1,q}(\Omega)$ with $q=\frac{p}{p-1}$, then up to a subsequence
$u_n \to u $ weakly in  $ W_0^{1,p}(\Omega)$ and $\tau_n \to \tau$  weakly in  $L^q(\Omega; \bR^N)$, 
where $u$ is a solution and $\tau$ is a momentum of the homogenized problem
\begin{eqnarray}\label{14}
\left\{
  \begin{array}{lll}
   \hfill  \tau(x) &\in\quad {\ss}_{hom}( \nabla u(x))& a.e. \,\,  x \in \Omega, \\
   -{\rm div} (\tau (x))&=\quad u^*(x) & a.e. \,\,  x \in \Omega.
  \end{array}
\right.
\end{eqnarray}
Here ${\ss}_{hom}$ can be defined variationally  as follows: for $\xi \in \bR^N$, ${\ss}_{hom}(\xi)=\partial \psi_{hom}(\xi)$, where
 \begin{equation}\label{15}
\psi_{hom}(\xi) := \min_{\phi \in W^{1,p}_\# (Q)}\frac{1}{|Q|}\int_Q \psi\big( x,\xi + \nabla \phi(x) \big) \dx,
\end{equation}
and 
\begin{equation}\label{16}
W^{1,p}_\# (Q)=\{u \in W^{1,p}(Q); \int_Q u(x) \, dx=0 \,\, \text{ and }\,\, u \text{ is } Q-\text{periodic}\}.
\end{equation}
As mentioned above, a similar result can be obtained for general maximal monotone maps  ${\ss}: \Omega \times \bR^N \to \bR^N$ with appropriate boundedness conditions (see below), by using the more cumbersome graph convergence (or $G$-convergence) methods. In this case, 
${\ss}_{hom}$ is defined by the following non-variational formula
\begin{equation}\label{beta.hom}
\hbox{${\ss}_{hom}(\xi)= \Big\{\int_{Q}g(y) \, dy \in \bR^N;\,  g\in L_\#^q(Q; \bR^N)$, $g(y)\in {\ss}(y, \xi+\nabla \psi(y))$  a.e. in $Q$ for some $\psi \in W^{1,p}_\#(Q)$\Big\}, }
 \end{equation}
where
\begin{equation}\label{16}
L_\#^q(Q; \bR^N):=\left\{g\in L^q(Q; \bR^N); \int_Q\langle g(y), \nabla \phi(y)\rangle_{\bR^N} \, dy=0  \text{ for every }  \phi \in W^{1,p}_\#(Q)\right\}.
\end{equation}
 More recently, the first-named author proposed a variational approach to deal with general maximal monotone operators, including corresponding equations of the form (\ref{eleven}) via the theory of selfdual Lagrangians on phase space \cite{Gh, G-book}. Our goal here is to describe how this approach is particularly well suited  to deal with the homogenization of such equations, first by showing that --just as in the case of a convex potential (\ref{twelve})-- the limiting process can be handled again through  $\Gamma$-convergence of associated selfdual Lagrangians, and secondly by giving a variational characterization for the limiting vector field (\ref{beta.hom}) in the same spirit as in (\ref{15}). 

 \quad We first recall that a {\it selfdual Lagrangian $L$} on a reflexive Banach space $X$, is any convex lower semi-continuous function on phase space $L:X\times X^*\to \R\cup\{+\infty\}$  that satisfy the following duality property:
  \begin{equation}
  L^*(u^*,u)=L(u, u^*) \quad \hbox{\rm for all $(u,u^*)\in X\times X^*$},
  \end{equation}
   where $X^*$ is the Banach space dual to $X$, and $L^*$ is the Fenchel-Legendre dual of $L$  in both variables, i.e.,
  \[
  L^*(u^*,u)= \sup  \{  \braket{ v}{ u^* }+   \braket{u}{v^* }-L(v, v^*): \,  (v,v^*)\in X\times X^*\}.
  \]
  Such Lagrangians satisfy the following basic property:
 \begin{equation*}\label{obs.1}
\hbox{$ L(u, u^*)-\langle u, u^*\rangle\geq 0$   for every $(u, u^*) \in X\times X^{*}$.}
 \end{equation*}
We then consider  the  corresponding  --possibly multivalued-- {\it selfdual vector field} $\bar \partial L: X\to 2^{X^*}$ defined
  for each $u \in X$ as the  --possibly empty-- subset  $\bar \partial L(u)$ of $X^*$ given by
 \begin{eqnarray}
\bar \partial L(u): = \{ u^* \in X^*; L(u,u^*)- \langle u,u^* \rangle=0 \}=\{ u^* \in X^*;  (u^*, u)\in \partial L(u,u^*) \}.
\end{eqnarray}
 Here $\partial L$ is the subdifferential of the convex function $L$ on $X\times X^*$, which should not be confused with $\bar \partial L$.

 Before going further, let us note that  {\it selfdual vector fields} are natural and far reaching extensions of subdifferentials of convex lower semi-continuous functions. Indeed, the most
basic selfdual Lagrangians are of the form $L(u,u^*)= \varphi (u)+\varphi^*(u^*)$ where  $\varphi$ is a convex  function in $X$, and $\varphi^*$ is its Fenchel dual  on $X^*$ (i.e., $
\phi^*(u^*)=\sup\{\langle u,u^* \rangle-\phi(u); \, u \in X\}$,  in which case  
\[
\bar \partial  L(u)= \partial \varphi (u).
\]
 More interesting  examples of selfdual Lagrangians are of the form $L(u,u^*)= \varphi (u)+\varphi^*(-\Gamma u+u^*)$  where $\varphi$  is a convex and lower semi-continuous
function on $X,$ and  $\Gamma: X\rightarrow X^*$ is a skew adjoint operator. The corresponding {\it selfdual vector field} is then
\[
\bar \partial L(u)=\Gamma u+ \partial \varphi (u).
\]  
Actually, both $\partial \varphi$ and $\partial \varphi +\Gamma$ are particular examples of the so-called {\it maximal monotone operators}, which are set-valued maps $\beta :X\to 2^{X^*}$ whose graph in $X\times X^*$ are maximal (for set inclusion) among all monotone subsets $G$ of $X\times X^*$, i.e., those $G$ satisfying
\begin{equation}
\hbox{$\langle x-y, p - q \rangle \geq 0$ for every $(x,p)$ and $(y, q)$ in $G$.}
\end{equation}
It turned out that the class of maximal monotone operators and the one of  selfdual vector fields coincide. Indeed, the following was proved in \cite{Gh}.
\begin{thm}\label{mm-sd} If ${\ss} : D({\ss})\subset X \to 2^{X^*}$ is a maximal monotone operator with a non-empty domain $D({\ss})$, then there exists a selfdual Lagrangian on $X \times X^*$ such that ${\ss} = \bar \partial L.$\\
Conversely, if $L$ is a proper selfdual Lagrangian on $X \times X^*,$ then the vector field $u \to \bar \partial L(u)$ is maximal monotone.
\end{thm}
 This  means that selfdual Lagrangians can be seen as the {\it potentials} for maximal monotone operators, in the same way as convex lower semi-continuous energies are the potentials of their own subdifferential, leading to a variational formulation and resolution of most equations involving maximal monotone operators such as the one in (\ref{eleven}). This was indeed done in \cite{Gh} in the case where $\beta$ does not depend on the state $x\in \Omega$. We shall however need to consider in this paper
 measurable families ${\ss}(x,.): \Omega \times \bR^N \to \bR^N$ of maximal monotone operators with suitable boundedness and coercivity conditions, and the possibility of associating to them measurable families $L(x, \cdot, \cdot):\Omega\times \bR^N \times \bR^N\to \bR^N$ of selfdual Lagrangians on $\bR^N \to \bR^N$ that reflect these conditions. For that we recall the definition of the class $M_{\Omega, p}(\bR^N)$ introduced in \cite{Ch-dalM-D}.
 
 \begin{defn}\label{M-omega} For a domain $\Omega$ in $\bR^N$, $p>1$ and $\frac{1}{p}+\frac{1}{q}=1$, we denote by $M_{\Omega, p}(\bR^N)$ the class of all possibly multi-valued functions ${\ss}: \Omega \times \bR^N \to \bR^N$ with closed values, which satisfy the following conditions:\\
(i) ${\ss}$ is measurable with respect to ${\cal L}(\Omega) \times {\cal B}(\bR^N)$ and ${\cal B}(\bR^N)$ where ${\cal L}(\Omega)$ is is the $\sigma$-field of all measurable subsets of $\Omega$ and ${\cal B}(\bR^N)$ is  the $\sigma$-field of all Borel subsets of $\bR^N.$\\
(ii) For a.e. $x \in \Omega$, the map ${\ss}(x,.): \bR^N \to \bR^N$ is maximal monotone. \\
(iii) There exist non-negative constants $m_1, m_2,c_1$ and $c_2$ such that for every $\xi \in \bR^N$ and $\eta  \in {\ss}(\xi)$,  
 
\begin{equation}
\langle \xi,\eta \rangle_{\bR^N}\geq \max \left\{\frac{c_1}{p}|\xi|^p-m_1, \frac{c_2}{q}|\eta|^q-m_2\right\}, 
\end{equation}
holds,  where $\langle .,. \rangle_{\bR^N}$ is the inner product in $\bR^N$.
\end{defn}
The following is the main application of the results in this paper.
\begin{thm}\label{main} Let $\Omega$ be a domain in $\bf R^N$, $q, p>1$ with $\frac{1}{p} + \frac{1}{q} = 1$, and assume $u^*_n \to u^* $ strongly in $W^{-1,q}(\Omega)$. Let  $u_n$ (resp., $\tau_n$) be (weak) solutions in $W_0^{1,p}(\Omega)$ (resp., momenta in $L^q(\Omega; \bR^N)$) for the Dirichlet boundary value problems (\ref{eleven}), 
where ${\ss}: \Omega \times \bR^N \to \bR^N$ belongs to $M_{\Omega, p}(\bR^N)$. 

If ${\ss}(., \xi)$ is  $Q$-periodic for an  open non-degenerate parallelogram $Q$ in  $\bR^n$  then, up to a subsequence
\[u_n \to u \qquad \text{ weakly in  }  W_0^{1,p}(\Omega), \]
\[\tau_n \to \tau  \qquad \text{ weakly in  } L^q(\Omega; \bR^N),\]
where $u$ is a solution and $\tau$ is a momentum of the homogenized problem
\begin{eqnarray}\label{mm-hom}
\left\{
  \begin{array}{lll}
    \tau(x) \in {\ss}_{hom}( \nabla u(x)) & a.e. \quad  x \in \Omega, \\
   -{\rm div} (\tau (x))= u^*(x) & a.e. \quad  x \in \Omega, \\
\, u \in W_0^{1,p}(\Omega).\\
  \end{array}
\right.
\end{eqnarray}
Here ${\ss}_{hom}=\bar \partial L_{hom}$, with $L_{hom}$ being a selfdual Lagrangian on $\bR^N \times \bR^N$ 
defined by
\begin{equation}\label{L.hom}
L_{hom}(a,b) := \min_{\substack{\phi \in W^{1,p}_\# (Q) \\ g \in L_\#^q(Q; \bR^N)}}\frac{1}{|Q|}\int_Q L\big( x,a + D\phi(x),b + g(x) \big) \dx,
\end{equation}
where for each $x\in \Omega$, 
$L(x, \cdot, \cdot)$ is a selfdual Lagrangian on $\bR^N \times \bR^N$ such that 
\begin{equation}\label{basic}
{\ss}(x, \cdot)=\bar \partial L (x, \cdot).
\end{equation}
\end{thm}
The above theorem will be a byproduct  of several results which have their own interest. In section 2, we consider various topologies on the class of selfdual Lagrangians that are relevant for homogenization. It turns out that the standard concept of $\Gamma$-convergence  is equivalent to the stronger notion of \emph{Mosco-convergence} in the context of selfdual Lagrangians. This has a direct implication on the corresponding maximal monotone operators. We also extend to selfdual Lagrangians one of the most attractive properties of the Mosco convergence of convex functions, which is that it implies  the  convergence of the graphs of their corresponding subdifferentials in the topology of Kuratowski-Painlev\'e on sets. We shall show in section 2 that similarly, the map $L \to \bar \partial L$ is continuous when we equip the class of selfdual Lagrangians with the topology of $\Gamma$-convergence and the class of maximal monotone operators with the topology of $G$-convergence.
 
 In section 3, we start by extending Theorem \ref{mm-sd} above by establishing a correspondence between state-dependent measurable maximal monotone operators in $M_{\Omega, p}(\bR^N)$ and  the following class of $\Omega$-dependent {\it selfdual} Lagrangian on $\Omega \times \mathbb{R}^{N} \times \mathbb{R}^{N}$. 

\begin{defn} \rm  An {\it ($\Omega, p)$-dependent {\it selfdual} Lagrangian on $\Omega \times \mathbb{R}^{N} \times \mathbb{R}^{N}$}  is a measurable 
 function $L: \Omega \times \mathbb{R}^{N} \times \mathbb{R}^{N}\to \R $
  such that 
 \begin{enumerate}
 \item   For any $x\in \Omega$, the map $(a,b)\to L(x, a,b)$ is a selfdual Lagrangian on $ \mathbb{R}^{N} \times \mathbb{R}^{N}.$
 \item There exist non-negative constants $C_0$ and $C_1$ and $n_0, n_1 \in L^1(\Omega)$ such that 
 \begin{equation}
 \hbox{$C_0(|a|^p + |b|^q - n_0(x)) \le L(x,a,b) \le C_1(|a|^p + |b|^q + n_1(x))$  for all $a,b \in \bR^N$.}
 \end{equation}

 \end{enumerate}
 \end{defn}
 
As in Theorem \ref{mm-sd}, any map ${\ss}: \Omega \times \bR^N \to \bR^N$ in $M_{\Omega, p}(\mathbb{R}^{N})$ can be seen as a potential of an  $\Omega$-dependent  selfdual Lagrangian  $L: \Omega \times \mathbb{R}^{N} \times \mathbb{R}^{N}\to \R,$ that is  $\bar \partial L(x,a)= {\ss}(x,a)$ for almost all $x \in \Omega.$

We then proceed to use the above representation of ${\ss}$ to give a variational resolution for the problem 
\begin{eqnarray}\label{div.form.0}
\left\{
\begin{array}{lll}
 f \in {\ss} (x, \nabla u(x))&  \text{ a.e. } x \in \Omega,\\
-{\rm div}(f) = u^*,\\
 u \in W^{1,p}_{0}(\Omega),
\end{array}
\right.
\end{eqnarray}
by ``lifting" the corresponding $\Omega$-dependent  selfdual Lagrangian $L$ on $\Omega \times \mathbb{R}^{N} \times \mathbb{R}^{N}$ to a selfdual Lagrangian on the function space $W_0^{1,p} (\Omega)\times W^{-1,q}(\Omega)$ via the formula:
\begin{equation}\label{lift00}
F(u,u^*) := \inf \{ \S L\big( x,\nabla u(x),f(x) \big) \, dx;\,  f \in L^q(\Omega; \bR^N), -{\rm div}(f) = u^*  \}.
\end{equation}
 A solution can then be obtained by simply minimizing for a given $u^* \in W^{-1,q}(\Omega)$ the non-negative functional 
 \[
I(u)=\inf_{ \substack{ f \in L^q(\Omega; \bR^N)\\ -{\rm div}(f) = u^* } }  \S \big [ L \big( x,\nabla u(x),f(x) \big)-\langle u(x), u^*(x) \rangle_{\bR^N} \big ] \, dx,
\]
 on $W^{1,p}_{0}(\Omega)$.  We end the section by showing that if  $\bar \partial L(x,.)= {\ss}(x,.)$,  then  
\begin{equation}\label{basic.2}
{\ss}_{hom}=\bar \partial L_{hom},
\end{equation}
 where ${\ss}_{hom}$ is defined in (\ref{beta.hom}) and $L_{hom}$ is as in  (\ref{L.hom}).
\\
We start section 4 by a homogenization result via $\Gamma$-convergence for general $Q$-periodic Lagrangians which are not necessarily selfdual. This is then applied to obtain the result claimed in Theorem \ref{main} above in the case of selfdual Lagrangians.  
    The last section is an appendix meant for auxiliary results that are needed throughout the paper.

\section{Preliminaries on selfdual Lagrangians}\label{Prelim}

We first recall the needed notions and results from the theory of selfdual Lagrangians developed in the book \cite{G-book}. We shall also establish new ones, in particular those regarding the convergence properties  in suitable topologies  of selfdual Lagrangians  and their associated maximal monotone vector fields. 
 $X$ will denote a real  reflexive Banach space and $X^*$ its dual.

\subsection{A variational principle for selfdual Lagrangians}

As mentioned in the introduction, maximal monotone operators ${\ss}$ can be written as ${\ss}=\bar \partial L$, 
where $L$ is a selfdual Lagrangian on $X\times X^*$, in such a way that solving the  equation 
 \begin{equation}
\hbox{$u^*\in  {\ss} (u)$,}
\end{equation}
amounts to minimizing the non-negative functional
$I(u):=L(u,u^*)-\langle u, u^*\rangle $. 
 The following  existence result is essential  for the sequel. It gives sufficient conditions for the infimum of selfdual Lagrangians  to be attained, and --as importantly-- to be zero.

\begin{thm}  \label{one} Let $L$ be a  selfdual Lagrangian  on a reflexive Banach space $X \times X^*$, let $u^*\in X^*$ be such that  $(0,u^*)\in Dom(L)$, and consider the functional $I (u):=L(u,u^*)-\langle u, u^*\rangle $.  Then  
  \begin{equation}
\inf_{u\in X}I (u)=0, 
\end{equation}
and in particular,  if the functional $I $ is coercive on $X$,  then  there exists $\bar u\in X$ such that   
\begin{equation}
\hbox{$
I (\bar u)=\min\limits_{u\in X}I(u)=0$ and 
$u^*\in \bar \partial L(\bar u)$.}
\end{equation}
  \end{thm}
Note that since $L_{u^*}(u, v^*):=L(u, u^*+v^*) -\langle u, u^*\rangle$ is a selfdual Lagrangian whenever $L$ is, it suffices to assume that $u^*=0$. The above theorem is then a consequence of  the following  result originally established in \cite{G1} (see also  \cite{G-book})  under a slightly stronger coercivity condition.

\begin{thm}\label{one0} Let $L$ be a  selfdual functional on a reflexive Banach space $X \times X^*$ such that  for some $u_{0}\in X$, the functional    $v^*\to L(u_{0},v^*)$ is bounded above on  a neighborhood of the origin in $X^{*}$.
  Then  there exists $\bar u\in X$ such that   $
I(\bar u)=\min\limits_{u\in X}I(u)=0.$
  \end{thm}
\textbf{Proof of Theorem \ref{one}}: Since $L$ is a selfdual Lagrangian on $X \times X^*$, so is its $\lambda$-regularization,
\begin{eqnarray*}
L_{\lambda}(u,u^*) =\inf \big \{   L(v,v^*)+ \frac{1}{2 \la}\|u-v\|^2+  \frac{\la}{2 }\|v\|^2+ \frac{1}{2 \la}\|u^*-v^*\|^2+  \frac{\la}{2 }\|v^*\|^2; \, v \in X, v^* \in X^*\big \},
\end{eqnarray*}
 for each $\lambda >0$,  by virtue of Lemma 3.2 in Chapter 2 of \cite{G-book}. Note that the Lagrangian $L_{\lambda}$ satisfies the boundedness condition of Theorem \ref{one0}. It then follows that
$\min_{u\in X}L_{\lambda}(u,0)=0.$ On the other hand, because of the properties of Yoshida regularization for convex functions, for each $(u,u^*) \in Dom (L)$ we have $\liminf_{\lambda\rightarrow 0}L_{\lambda}(u,u^*)=L(u,u^*).$ It follows that
\[\inf_{u\in X}L(u,0)=\inf_{u\in X}\liminf_{\lambda\rightarrow 0}L_{\lambda}(u,0)=\liminf_{\lambda\rightarrow 0}\inf_{u\in X}L_{\lambda}(u,0)=0.\]
Therefore $\inf_{u\in X}I(u)=0$. Now if $I$ is coercive then the minimum is attained for some $\bar u \in X,$ i.e.,   $I(\bar u)=L(\bar u,0)=0$ and consequently 
$\bar u$ is  a solution of 
$0\in \bar \partial L(\bar u)$.
\hfill $\square$

\subsection{Mosco and $\Gamma$-convergence of selfdual functionals}
 We  first recall  the main definitions and statements in the theory of  variational convergence for functionals, as well as the graph convergence for possibly multi-valued operators. 
 A complete study relating the various modes of convergence of convex functions and their subdifferentials can be found in \cite{D-M-S-Mosco}.  

\begin{defn}  Let $F_n$ and $F$ be functionals on a reflexive Banach space $X$.
The sequence $\{ F_n \}$ is said to \emph{$\Gamma$-converge} (resp., \emph{Mosco-converge}) to $F$,  
if the following two conditions are satisfied:
\begin{enumerate}
\item For any sequence $\{ u_n \} \subset X$ such that $u_n \to u$ strongly (resp., $u_n \rightharpoonup u$ weakly) in $X$ to some $u \in X$, one has
\[
  F(u) \le \liminf_{n\to \infty}F_n(u_n).
\]
\item For any $u \in X$, there exists a sequence $\{ u_n \} \subset X$ such that $u_n \to u$ strongly in $X$ and
\[
\lim_{n\to \infty} F_n(u_n) = F(u).
\]
\end{enumerate}
\end{defn}

 The following is a fundamental property of Mosco-convergence.
\begin{lem}\label{CMC}
Let $F_n, F$ be a proper convex lower semi-continuous functionals, then $\{ F_n \}$ Mosco-converge to $F$ if and only their Fenchel-Legendre duals $\{ F^*_n \}$ Mosco-converge to $F^*$. 
\end{lem}
In the following we note that this property implies the agreable fact that Mosco and $\Gamma$-convergence are actually equivalent for a sequence of selfdual Lagrangians $\{L_n\}$,  as long as the limiting Lagrangian $L$ is itself selfdual.

 \begin{thm}\label{sdmosdef} Let $\{ L_n \}$ be a family of selfdual Lagrangians on $X\times X^*$, where $X$ is a reflexive Banach space, and let $L$ be a Lagrangian on $X\times X^*$. The following statements are then  equivalent:
 \begin{enumerate}
\item $\{ L_n \}$ Mosco-converges to $L$.
\item  $L$ is selfdual and $\{ L_n \}$ $\Gamma$-converges to $F$.
 \item  $L$ is selfdual and for any $(u,u^*) \in X\times X^*$, there exists a sequence $(u_n,u^*_n)$ converging strongly to $ (u,u^*)$ in $X\times X^*$ such that
 \[
 \limsup_n L_n(u_n,u^*_n) \le L(u,u^*).
 \]
 \end{enumerate}
 \end{thm}

 \textbf{Proof.} For $(1) \rightarrow (2)$ we just need to prove  that $L$ is selfdual since Mosco convergence clearly implies $\Gamma$-convergence.  Since  $L$ is the Mosco limit of $ L_n$, it follows from  Lemma  \ref{CMC} that  $L^*$ is a Mosco limit of $ L^*_n$. Denoting
 \[
\hbox{$ L_n^T (u^*,u) := L_n(u,u^*)$  and  $ L^T (u^*,u) := L(u,u^*)$,}
 \]
it follows that  $L^T$ is a Mosco-limit of $L^T_n$ on $X^*\times X$. On the other hand, by selfduality of $L_n$ we have  that $L_n^T = L_n^*$ from which we obtain  that $L^T =\lim_n L^T_n =\lim_n L^*_n = L^*$, and therefore
$
L^T = L^*,
$
 and $L$ is therefore selfdual.\\
(2)$\to$(3) follows from the definition of $\Gamma$-convergence.\\
For (3)$\to$(1) we let  $(u^*,u) \in X^* \times X$ and consider a sequence $\{ (u^*_n,u_n) \} \subset X^*\times X$ such that $(u^*_n,u_n) \rightharpoonup (u^*,u)$ weakly in $X^* \times X$. By the definition of Fenchel-Legendre duality we have
\begin{equation}\label{mo200}
\liminf_n L^*_n(u^*_n,u_n) = \liminf_n \sup_{(v,v^*)\in X\times X^*} \{ \langle u_n,v^* \rangle + \langle v,u^*_n \rangle - L_n(v,v^*) \}.
\end{equation}
Consider now an arbitrary pair $(\tilde{u}, \tilde{u}^*)$ and let $\{ (\tilde{u}_n, \tilde{u}^*_n) \}$ be the recovery sequence given in item (3). It follows from (\ref{mo200}) that
\[
\liminf_n L_n(u_n,u^*_n) \ge \liminf_n \Big( \langle u_n,\tilde{u}^*_n \rangle + \langle \tilde{u}_n,u^*_n \rangle - L_n(\tilde{u}_n,\tilde{u}^*_n) \Big) = \langle u,\tilde{u}^* \rangle + \langle \tilde{u},u^* \rangle - L(\tilde{u},\tilde{u}^*).
\]
Since $(\tilde{u}, \tilde{u}^*)$ is arbitrary, taking the supremum over all  $(\tilde{u}, \tilde{u}^*)$ yields
\[
\liminf_n L^*_n(u^*_n,u_n) \ge L^*(u^*,u).
\]
  Since both $L_n$ and $L$ are selfdual, this implies that
\[
\liminf_n L_n(u_n, u^*_n) \ge L(u, u^*),
\]
 and therefore
that $L$ is a Mosco-limit of $L_n$.
\hfill $\square$

\begin{remark} Note that while the Mosco convergence of selfdual Lagrangians automatically implies that  the  limiting Lagrangian $L$ is itself selfdual,  this fails  for $\Gamma$-convergence as shown in the following example.
\end{remark}

 Let  $H$ be an infinite dimensional Hilbert space. Consider a set $\{e_n\}$ with $\|e_n\| = 1$ and $e_n \rightharpoonup 0$ (For example, the orthonormal basis of the space). Define
\[
L_n(u,u^*) := \frac{1}{2}\|u-e_n\|^2 + \frac{1}{2}\|u^*\|^2 + \langle u^*,e_n \rangle, 
\]
 in such a way that $L_n$ is selfdual. It can be checked directly that for any $(u_n,u^*_n) \to (u,u^*)$ in $H \times H$ we have $
 \lim_n L_n(u_n,u^*_n) = L(u,u^*)$, 
where
 \[
L(u,u^*) := \frac{1}{2}\|u\|^2 + \frac{1}{2}\|u^*\|^2 + \frac{1}{2}.
 \]
  This means that $L$ is a $\Gamma$-limit  of $L_n$. On the other hand, it is easily seen that  $L$ is not selfdual and therefore we do not have Mosco convergence.

\subsection{Continuity of $L\to \bar \partial L$ for the $\Gamma$-convergence of selfdual Lagrangians} 

One of the most attractive properties of Mosco convergence is the fact that for convex functions it implies the graph convergence (or \emph{Kuratowski-Painlev\'e  convergence}) of their corresponding subdifferentials \cite[Theorem 4.2]{Att-B}. We shall extend this result to selfdual Lagrangians by showing that their Mosco (or $\Gamma$-convergence) also yield the graph convergence of their derived vector fields (i.e., their corresponding maximal monotone operators). 

 Considering a sequence of sets $\{ A_n \}$ in $X$, the corresponding sequential lower and upper limit sets are respectively given by
\[
Li_X \big( A_n \big) = \{ u \in X \; : \; \exists u_n \to u,\, u_n \in A_n \}, 
\]
and
\[
Ls_X \big( A_n \big) = \{ u \in X \; : \; \exists k(n) \to \infty, \, \exists u_{n(k)} \to u,\, u_{n(k)} \in A_k \}.
\]
In other words,  $Li \big( A_n \big)$ corresponds to the collection of all \emph{limit} points of the sequence $\{ A_n \}$ and $Ls \big( A_n \big)$ is the collection of all \emph{cluster} points of the sequence $\{ A_n \}$.
We clearly have $
  Li_X (A_n) \subseteq Ls_X (A_n)$.

\begin{defn}
 A sequence of subsets $\{ A_n \}$ of $X$ is said to converge to $A \subset X$, in the sense of \emph{Kuratowski-Painlev\'e}, if
$Ls_X (A_n) = A = Li_X (A_n).$
\end{defn}

 This definition, when $X$ is replaced by the phase space  $X\times X^*$ and when the subsets $A_n$ are graphs of maps from $X$ to $X^*$, is also refered to as \emph{graph}-convergence  (see Definition 3.5 on \cite{Ch-dalM-D}). 
 
Recall that for a selfdual Lagrangian $F$ on $X \times X^*,$ its associated vector field at $u \in X$ is denoted by $\bar \partial F (u)$ and  given by
$
\bar \partial F (u)=\{u^* \in X^*; F(u,u^*)=\langle u, u^*\rangle\}$. We shall therefore  also denote by $ \bar \partial F$ the graph of $ \bar \partial F$ in $ X \times X^*$,  i.e.,
\[
\hbox{$(u,u^*) \in \bar \partial F$ \text{ if and only if } $u^* \in \bar \partial F(u)$.}
\]
 The following is the main result of this section.
 \begin{thm}\label{KP} Let $X$ be a reflexive Banach space and suppose $\{F_n\}$ is a family of selfdual Lagrangians on $X \times X^*.$
If $F:X \times X^* \to \bR \cup \{+\infty\}$ is a  selfdual Lagrangian that is a $\Gamma$-limit of $\{F_n\}$, then  
the graph of $ \bar \partial F_n $ converge to the graph of $ \bar \partial F$ in the sense of Kuratowski-Painlev\'e.
\end{thm}

For the proof, we shall make use of the following theorem that can be seen as the counterpart of the  \emph{Br\o ndsted-Rockafellar} result for convex functions \cite{phelps}.

\begin{lem}\label{BRL}
Let $L: X \times X^* \to \bR\cup \{+\infty\}$ be a selfdual Lagrangian and 
assume that for a pair $(u_0,u^*_0) \in X\times X^*$, we have 
$
L(u_0,u^*_0) - \langle u_0,u^*_0 \rangle \le \e.
$
Then, there exists a pair $(u_\e,u^*_\e) \in \dL$ such that
\begin{enumerate}
\item $\| u_\e - u_0 \| \le \sqrt{\e},$
\item $\| u^*_\e - u^*_0 \|_* \le \sqrt{\e},$
\item $| L(u_\e,u^*_\e) - L(u_0,u^*_0) | \le 2\e + \sqrt{\e}(\|u_0\| + \|u^*_0\|_*)$.
\end{enumerate}
\end{lem}
\textbf{Proof:}  First assume  that $M$ is a selfdual Lagrangian such that $
M(0,0) \le \e$. We claim that there exists then a pair $(v_\e,v^*_\e) \in \dM$ such that
\begin{enumerate}
\item $\| v_\e \| \le \sqrt{\e},$
\item $\| v^*_\e \|_* \le \sqrt{\e}$,
\item $|M(v_\e,v^*_\e)| \leq \e$.
\end{enumerate}
Indeed,  denote by $J$ the \emph{duality mapping} from $X$ to $X^*$ and use the fact that 
$\dM$ is a maximal monotone operator to find  $\tilde{u} \in X$ such that
\[
-J\tilde{u} \in \dM(\tilde{u}).
\]
It follows that $M(\tilde{u},-J\tilde{u}) = \langle \tilde{u},-J\tilde{u} \rangle = -\| \tilde{u} \|^2$.
Now, since $M$ is selfdual, we have
\[
M(0,0) =M^*(0,0)= \sup_{(u,u^*) \in X\times X^*} -M(u,u^*) \ge -M(\tilde{u},-J\tilde{u}) = \| \tilde{u} \|^2,
\]
from which we obtain that 
$
\| \tilde{u} \|^2 \le \e.
$
Since $\| \tilde{u} \| = \| J\tilde{u} \|_* $,  it suffices to set 
$
v_\e := \tilde{u}
$
and
$
v^*_\e := J\tilde{u}$, to obtain that $\| v_\e \| = \| v^*_\e \|_*\le \sqrt{\e}$ and $|M(v_\e,v^*_\e)|= \| v_\e \|^2\le \e.$

To complete the proof of Theorem \ref{BRL},  we  set
\[
M(u,u^*) := L(u + u_0, u^* + u^*_0) - \langle u , u^*_0 \rangle - \langle u_0 , u^* \rangle - \langle u_0 , u^*_0 \rangle,
\]
which is a selfdual Lagrangian on $X \times X^*$. The hypothesis yields that
\[
M(0,0) = L(u_0,u^*_0) - \langle u_0,u^*_0 \rangle \le \e.
\]
It then follows from the above that there exists  a pair $(v_\e,v^*_\e)\in \dM$ such that
 $\| v_\e \| \le \sqrt{\e},$
 $\| v^*_\e \|_* \le \sqrt{\e}$ and
 $|M(v_\e,v^*_\e)| \leq \e.$ Setting
$u_\e := v_\e + u_0$ and
$
u^*_\e := v^*_\e + u^*_0, 
$
and since  $M(v_\e,v^*_\e) = \langle v_\e , v^*_\e \rangle$, we have
$
L(u_\e,u^*_\e) = \langle u_\e , u^*_\e \rangle,
$
and therefore $(u_\e,u^*_\e) \in \dL.$ Note also that $\| u_\e - u_0\| \le \sqrt{\e}$ and $\| u^*_\e - u^*_0\|_* \le \sqrt{\e}$.
Finally, we have
\[
L(u_\e,u^*_\e) - L(u_0,u^*_0) = M(v_\e,v^*_\e) +\langle v_\e,u^*_0 \rangle + \langle u_0, v^*_\e \rangle - M(0,0), 
\]
which together with   $| M(v_\e,v^*_\e) | \le \e$, yields that
\[
| L(u_\e,u^*_\e) - L(u_0,u^*_0) | \le 2\e + \sqrt{\e}(\|u_0\| + \|u^*_0\|).
\]
\hfill $\square$

\textbf{Proof of Theorem \ref{KP}.}  Fix $(u,u^*) \in \dF$. There exists then in view of the $\Gamma$-convergence, a sequence $(u_n,u^*_n)$ converging strongly to $(u, u^*)$ in $ X \times X^*$ such that $F_n(u_n,u^*_n) \to F(u,u^*)$.
 We  then have
$
F(u,u^*) = \langle u,u^* \rangle = \lim_n \langle u_n,u^*_n \rangle,
$ 
and therefore  if we define $
\e_n := F_n(u_n,u^*_n) - \langle u_n,u^*_n \rangle,
$
we obtain that $
\lim_n \e_n = 0$.  Hence, by Lemma  \ref{BRL}, we have the existence of a pair $(\tilde{u}_n,\tilde{u}^*_n) \in \dF_n$ such that $\| u_n - \tilde{u}_n \| < \sqrt{\e_n}$ and $\| u^*_n - \tilde{u}^*_n \|_* < \sqrt{\e_n}$. Clearly $\tilde{u}_n \to u$ and $\tilde{u}^*_n \to u^*$ as $\e_n \to 0$.
This shows that $\dF \subset Li(\dF_n)$.

 To complete the proof, we just need to show that $Ls(\dF_n) \subset \dF$. Letting
 $(v,v^*) \in Ls(\dF_n)$, there exists some sequence $(v_{n_k},v^*_{n_k}) \in \dF_{n_k}$ such that  $(v_{n(k)},u_{n(k)}) \to (v,v^*)$. Now take an arbitrary $(u,u^*) \in \dF$. From what we have shown, there exists a sequence $(u_n,u^*_n) \in \dF_n$ such that $(u_n,u^*_n) \to (u,u^*)$. For each $k$ we have
$
\langle u_{n(k)} - v_{n(k)} , u^*_{n(k)} - v^*_{n(k)} \rangle \ge 0,
$
and as $k \to \infty$ we get
$
\langle u - v, u^* - v^* \rangle \ge 0.
$
The above holds  for all $(u,u^*) \in \dF$ and so by the maximality of $\dF$ we obtain that $(v,v^*) \in \dF$, which completes the proof.
\hfill $\square$

\section{A selfdual variational approach to existence theory}

In  this section,  we first establish a correspondence between maximal monotone maps in $M_{\Omega, p}(\bR^N)$ and a class of  $\Omega$-dependent selfdual Lagrangians. We then proceed to give a variational formulation and resolution  to  equation (\ref{eleven}) even in the case where the maximal monotone operator $\beta$ is nor derived from the potential of a convex function.  

\subsection{Selfdual Lagrangians associated to maximal monotone operators}

\begin{defn} \rm Let $\Omega$ be a domain in  $\mathbb{R}^{N}$. 

\quad (i) A   function $L: \Omega \times \mathbb{R}^{N} \times \mathbb{R}^{N}\to \R \cup \{+\infty\}$ is said to be an  {\it $\Omega$-dependent Lagrangian} on $\Omega \times \mathbb{R}^{N} \times \mathbb{R}^{N}$,  if it is  measurable with respect to   the  $\sigma$-field  generated by the products of Lebesgue sets in  $\Omega$ and Borel sets in $\mathbb{R}^{N} \times \mathbb{R}^{N}$.
   
\quad (ii)   Such a Lagrangian $L$ is  said to be {\it selfdual  on $\Omega \times \mathbb{R}^{N} \times \mathbb{R}^{N}$}  if  for any $x\in \Omega$, the map $L_x:(a,b)\to L(x, a,b)$ is a selfdual Lagrangian on $ \mathbb{R}^{N} \times \mathbb{R}^{N},$
i.e., if $L^*(x,b,a)=L(x,a,b)$ for all $a,b \in \bR^N$ where
\[
  L^*(x,b,a)= \sup  \{ \langle b, \xi \rangle_{\bR^N}+ \langle a, \eta \rangle_{\bR^N}-L(x,\xi, \eta): \,  (\xi, \eta)\in \mathbb{R}^{N} \times \mathbb{R}^{N}\}.
  \]
\end{defn}
The following was proved in \cite{Gh} for a single maximal monotone operator.

\begin{pro} \label{rep} If ${\ss} \in M_{\Omega, p}(\bR^N)$ for some $p>1$, then there exists an $\Omega$-dependent selfdual Lagrangian $L:\Omega \times \bR^N \times \bR^N \to \bR$ such that ${\ss}(x, .)=\bar \partial L(x,.)$ for a.e. $x \in \Omega$ and
\begin{equation}\label{est200}
C_0(|a|^p + |b|^q - n_0(x)) \le L(x,a,b) \le C_1(|a|^p + |b|^q + n_1(x)) \qquad \text{ for all } a,b \in \bR^N.
\end{equation}
where $C_0$ and $C_1$ are two positive constants and $n_0, n_1 \in L^1(\Omega).$\\
Conversely, if $L:\Omega \times \bR^N \times \bR^N \to \bR$ is an $\Omega$-dependent selfdual Lagrangian satisfying (\ref{est200}), then $\bar \partial L(x,.) \in M_{\Omega, p}(\bR^N).$
\end{pro}
\textbf{Proof.}  Let $N:\Omega \times \bR^N \times \bR^N \to \bR \cup\{+\infty\}$ be the Fitzpatrick function \cite{fitz} associated to ${\ss}$, i.e.,
\[N(x,a,b):=\sup\{\langle b, \xi \rangle_{\bR^N}+\langle a-\xi, \eta \rangle_{\bR^N}; \eta \in {\ss}(x, \xi)\}.\]
Note that measurability assumptions on ${\ss}$ ensures that $N$ is a normal integrand. Also, by the properties of  the Fitzpatrick function \cite{G-book}, it follows that
\[
\hbox{$N^*(x,b,a)\geq N(x,a,b) \geq \langle a, b\rangle_{\bR^N}$ for a.e. $x \in \Omega$ and for all $a,b \in \bR^N$.}
 \]
Moreover, 
\begin{equation}\label{fit201}\
\hbox{$\eta \in  {\ss}(x, \xi)$ if and only if  $N^*(x,\eta ,\xi)=N(x,\xi,\eta)=\langle \eta, \xi \rangle_{\bR^N}$  a.e. $x \in \Omega$.}
\end{equation}
Define $L:\Omega \times \bR^N \times \bR^N \to \bR$ by
\[L(x,a,b)=\inf \big \{  \frac{1}{2} N(x, a_1, b_1)+ \frac{1}{2} N^*(x, b_2, a_2)+\frac{1}{4p}|a_1-a_2|^p+\frac{1}{4q}|b_1-b_2|^q; (a,b)=\frac{1}{2}(a_1,b_1)+\frac{1}{2}(a_2,b_2) \big\}.\]
We shall show that $L$ is   $\Omega$-dependent selfdual Lagrangian such that
\begin{equation}\label{fit200}
\hbox{$N^*(x,b,a)\geq L(x, a,b) \geq N(x,a,b)$  for a.e. $x \in \Omega$  and for all  $a,b \in \bR^N$.}
\end{equation}
Fix $a,b \in \bR^N.$ We have
\begin{eqnarray*}
L^*(x,b,a)&=& \sup_{\xi, \eta \in \bR^N}\{\langle \xi, b\rangle_{\bR^N}+\langle a, \eta\rangle_{\bR^N}-L(x, \xi, \eta) \}\\
&=& \sup_{\xi, \eta \in \bR^N}\Big \{\langle \xi, b\rangle_{\bR^N}+\langle a, \eta\rangle_{\bR^N}-\frac{1}{2} N(x, \xi_1, \eta_1)- \frac{1}{2} N^*(x, \xi_2, \eta_2)\\&&-\frac{1}{4p}|\xi_1-\xi_2|^p-\frac{1}{4q}|\eta_1-\eta_2|^q; \quad(\xi,\eta)=\frac{1}{2}(\xi_1,\eta_1)+\frac{1}{2}(\xi_2,\eta_2) \Big\}\\
&=& \frac{1}{2}\sup_{\xi_1,\xi_2, \eta_1, \eta_2 \in \bR^N}\Big \{\langle \xi_1+\xi_2, b\rangle_{\bR^N}+\langle a, \eta_1+\eta_2\rangle_{\bR^N}- N(x, \xi_1, \eta_1)-  N^*(x, \xi_2, \eta_2)\\&&-\frac{1}{2p}|\xi_1-\xi_2|^p-\frac{1}{2q}|\eta_1-\eta_2|^q \Big\}.
\end{eqnarray*}
Using the fact that the Fenchel dual of some of two functions is their inf-convolution, we obtain
\begin{eqnarray*}
L^*(x,b,a)= \frac{1}{2}\inf_{a_1,b_1 \in \bR^N}\Big \{ N^*(x, b_1, a_1)+ N(x, 2a-a_1, 2b-b_1)+\frac{2^{q-1}}{q}|b-b_1|^q+\frac{2^{p-1}}{2p}|a-a_1|^p \Big\}.
\end{eqnarray*}
Setting $a_2=2a-a_1$ and $b_2=2b-b_1$ we have $a=\frac{a_1+a_2}{2}$ and $b=\frac{b_1+b_2}{2}.$ It then follows that
\begin{eqnarray*}
L^*(x,b,a)&=& \frac{1}{2}\inf_{a_1,b_1, a_2, b_2 \in \bR^N}\Big \{ N^*(x, b_1, a_1)+ N(x, a_2, b_2)\\&&+\frac{2^{q-1}}{q}|\frac{b_1-b_2}{2}|^q+\frac{2^{p-1}}{2p}|\frac{a_1-a_2}{2}|^p; \quad  (a,b)=\frac{1}{2}(a_1,b_1)+\frac{1}{2}(a_2,b_2)\Big\}\\
&=& \inf \Big \{ \frac{1}{2} N^*(x, b_1, a_1)+ \frac{1}{2} N(x, a_2, b_2)+\frac{1}{4q}|b_1-b_2|^q+\frac{1}{4p}|a_1-a_2|^p;\\&& \quad  (a,b)=\frac{1}{2}(a_1,b_1)+\frac{1}{2}(a_2,b_2)\Big\}\\
&=& L(x,a,b).
\end{eqnarray*}
Thus, $L$ is a   $\Omega$-dependent selfdual Lagrangian. Inequalities (\ref{fit200}) simply follow from the definition and selfduality of  $L$. We shall now prove that $L$ satisfies the estimate (\ref{est200}). Note first that for all  $\eta \in {\ss}(x, \xi)$
we have
\[\frac{1}{p}|\xi|^p+\frac{1}{q}|\eta|^p\leq m_1+m_2+ (c_1+c_2)\langle \xi, \eta \rangle_{\bR^N}.\]
 It follows from  the definition of the Fitzpatrick function $N$ that
\begin{eqnarray}\label{fit202}
N(x,a,b)&=&\sup\{\langle b, \xi \rangle_{\bR^N}+\langle a-\xi, \eta \rangle_{\bR^N}; \eta \in {\ss}(x, \xi)\}\nonumber\\
&\leq & \sup \Big \{\langle b, \xi \rangle_{\bR^N}+\langle a, \eta \rangle_{\bR^N}-\frac{1}{p(c_1+c_2)}|\xi|^p-\frac{1}{q(c_1+c_2)}|\eta|^q-\frac{m_1+m_2}{c_1+c_2}; \eta \in {\ss}(x, \xi)\Big \} \nonumber\\
&\leq & \sup_{\xi, \eta \in \bR^N}\Big \{\langle b, \xi \rangle_{\bR^N}+\langle a, \eta \rangle_{\bR^N}-\frac{1}{p(c_1+c_2)}|\xi|^p-\frac{1}{q(c_1+c_2)}|\eta|^q-\frac{m_1+m_2}{c_1+c_2}\Big \}\nonumber \\
&=&\frac{(c_1+c_2)^{p-1}}{p}|a|^p+\frac{(c_1+c_2)^{q-1}}{q}|b|^q+\frac{m_1+m_2}{c_1+c_2}.
\end{eqnarray}
Let $\eta_0(x) \in {\ss}(x,0).$ By assumption $|\eta_0(x)|^q\leq m_2+ \langle 0,\eta_0(x) \rangle=m_2$ for a.e. $x \in \Omega,$ from which we get $\eta_0 \in L^q(\Omega).$ It also follows from (\ref{fit201}) that $N^*(x, \eta_0(x),0)=0$ for a.e. $x \in \Omega.$  From the definition of $L$ and (\ref{fit202}), we get that
\begin{eqnarray*}
L(x,a,b) &\leq& \frac{1}{2}N(x, 2a-\eta_0(x), 2b)+ \frac{1}{2} N^*(x, \eta_0(x),0) +\frac{2^q}{4q}|b|^q+\frac{2^p}{4p}|a-\eta_0(x)|^p\\
&\leq& C_1 (|a|^p+|b|^q+n_1(x)) \,\, a.e. \quad  x \in \Omega,
\end{eqnarray*}
where $C_1$ is a positive constant and $n_1 \in L^1(\Omega).$  The reverse inequality follows from the selfduality of $L.$\\

Conversely, let $L$ be a $\Omega$-dependent selfdual Lagrangian satisfying  (\ref{est200}). If $\eta \in \bar \partial L(x, \xi)$ then
\[\langle \xi, \eta \rangle= L(x, \xi, \eta)\geq C_0(|\xi|^p+|\eta|^q-n_0(x)),\]
from which we conclude that $\bar \partial L(x,.) \in M_{\Omega, p}(\bR^N).$ \hfill $\square$

\subsection{Self-dual Lagrangians on $W^{1,p}_{0}(\Omega) \times W^{-1,q}(\Omega)$} \label{SDLags}

We now show how one can  ``lift" an $\Omega$-dependent selfdual Lagrangian  to a selfdual Lagrangian on the phase space  $W^{1,p}_{0}(\Omega) \times W^{-1,q}(\Omega).$ 
This will allow us to give a variational formulation and resolution --via Theorem \ref{one}--  of equations involving  maximal monotone operators in divergence form. The following  extends a result  in \cite{Gh}.

\begin{thm} Let  ${\ss} \in M_{\Omega, p}(\bR^N)$ for some $p>1$, then for every $u^* \in W^{-1,q}(\Omega)$ with $\frac{1}{p}+ \frac{1}{q}=1$, there exist  $\bar u \in W^{1,p}_{0}(\Omega)$ and $\bar f(x) \in L^q(\Omega; \bR^N)$ such that
\begin{eqnarray}\label{div.form}
\left\{
\begin{array}{lll}
\bar f \in {\ss} (x, \nabla \bar u(x))&  \text{ a.e. } x \in \Omega\\
-{\rm div}(\bar f) = u^*.
\end{array}
\right.
\end{eqnarray}
 It is obtained by minimizing the functional 
 \[I(u):=\inf_{ \substack{ f \in L^q(\Omega; \bR^N) \\ -{\rm div}(f) = u^* } }\S \big [ L \big( x,\nabla u(x),f(x) \big) -\langle u(x), u^*(x) \rangle_{\bR^N} \big ]\, dx
 \]
 on  $W^{1,p}(\Omega)$, where $L$ is an $\Omega$-dependent selfdual Lagrangian on $\Omega \times \mathbb{R}^{N} \times \mathbb{R}^{N}$ associated to ${\ss}$ in such a way that $\bar \partial L (x, \cdot)={\ss}(x, \cdot)$ for a.e $x\in \Omega$.

\end{thm}

The above theorem will follow from the representation of a maximal monotone map in $M_{\Omega, p}(\bR^N)$ by an $\Omega$-dependent selfdual Lagrangian on $\Omega \times \mathbb{R}^{N} \times \mathbb{R}^{N}$ (Proposition \ref{rep}) combined with the following two propositions. 
\begin{pro}\label{lift}
Suppose $L$ is an $\Omega$-dependent selfdual Lagrangian on $\Omega \times \mathbb{R}^{N} \times \mathbb{R}^{N}$ such that $L(\cdot ,0,0) \in L^1(\Omega)$, then the Lagrangian defined on $W^{1,p}_{0}(\Omega) \times W^{-1,q}(\Omega)$ by
\begin{eqnarray}\label{lift0}
F(u,u^*) := \inf \{ \S L\big( x,\nabla u(x),f(x) \big) \, dx; f \in L^q(\Omega; \bR^N), -{\rm div}(f) = u^*  \},
\end{eqnarray}
is selfdual.

\end{pro}
\textbf{Proof:}
 Denote $W^{1,p}_{0}(\Omega)$ by $X$ and its dual $W^{-1,q}(\Omega)$ by $X^*$.  For a fixed $(v^*,v) \in X^* \times X $,  we have
\begin{eqnarray*}
F^*(v^*,v) &=& \sup \{ \langle u,v^* \rangle + \langle u^*,v \rangle - F(u,u^*);u \in X,  u^* \in X^* \}\\
&=& \sup_{ \substack{ f \in L^q(\Omega; \bR^N) \\ -{\rm div}(f) = u^* \\ (u,u^*) \in X\times X^*}}\left\{ \langle u,v^* \rangle + \langle u^*,v \rangle - \S L\big( x,\nabla u(x),f(x) \big) \dx\right\} \\
&=& \sup_{ \substack{  }} \{ \langle u,v^* \rangle - \langle {\rm div}(f),v \rangle - \S L\big( x,\nabla u(x),f(x) \big) \dx;u \in X, f \in L^q(\Omega; \bR^N) \}\\
&=& \sup \{ \langle u,v^* \rangle + \langle f,\nabla v \rangle - \S L\big( x,\nabla u(x),f(x) \big) \dx; u \in X, f \in L^q(\Omega; \bR^N)\}.
\end{eqnarray*}
Now set
$
E:= \{ g \in L^p(\Omega; \bR^N); g = \nabla u, \, u \in X \}
$
and let $\chi_E$ be the indicator function in $L^p(\Omega; \bR^N)$, e.g.,
\[
\chi_E (g) = \left\{ \begin{array}{lc} 0 & g \in E, \\ +\infty & \mbox{elsewhere.} \end{array}\right.
\]
An easy computation shows that
\[
\chi^*_E (f) = \left\{ \begin{array}{lc} 0 & {\rm div}(f) = 0, \\ +\infty & \mbox{elsewhere.} \end{array}\right.
\]
Fix   $f_0 \in L^q(\Omega; \bR^N)$ with  $-{\rm div}(f_0) = v^*$.
It follows that
\begin{eqnarray*}
F^*(v^*,v)&=& \sup \{ \langle g,f_0 \rangle + \langle f,\nabla v \rangle - \S L\big( x,g(x),f(x) \big) \dx - \chi_E (g); g \in L^p(\Omega; \bR^N), f \in L^q(\Omega; \bR^N)  \} \\
&=& \inf \{ \S L^*(x,f_0 - f, \nabla v) \dx + \chi^*_E (f); f \in L^q(\Omega; \bR^N) \}.
\end{eqnarray*}
Note that we have used the fact that $\big ( \S L\big( x,.,. \big) \dx \big )^*(g,f)= \S L^*\big( x,f(x),g(x) \big) \dx$ that holds since $L(.,0,0) \in L^1(\Omega).$
We finally get
\begin{eqnarray*}
F^*(v^*,v) &=& \inf\{ \S L^*(x,f_0 - f, \nabla v) \dx; f \in L^q(\Omega; \bR^N),{\rm div}(f) = 0 \}\\
&=& \inf\{ \S L(x, \nabla v, f_0 - f) \dx; f \in L^q(\Omega; \bR^N), {\rm div}(f) = 0 \} \\
&=&\inf \{ \S L(x, \nabla v, f) \dx; f \in L^q(\Omega; \bR^N),-{\rm div}(f) = v^* \}\\
& =& F(v,v^*).
\end{eqnarray*}
\hfill $\square$

Here is our variational resolution for equation ({\ref{div.form}).

\begin{pro} \label{prop.1} Suppose $L$ is  $\Omega$-dependent selfdual Lagrangian on $\Omega \times \mathbb{R}^{N} \times \mathbb{R}^{N}$
Assume the following coercivity  condition:
 \begin{equation}\label{A}
\hbox{$L(x,a,b)\geq m(x)+ C(|a|^p+|b|^q)$ for all $a,b \in \bR^N$,}
\end{equation}
where $m \in L^1(\Omega)$ and $C$ is a positive constant.
Then for every $u^* \in W^{-1,q}(\Omega)$ the functional 
\[
I(u)=\inf_{ \substack{ f \in L^q(\Omega; \bR^N)\\ -{\rm div}(f) = u^* } }  \S \big [ L \big( x,\nabla u(x),f(x) \big)-\langle u(x), u^*(x) \rangle_{\bR^N} \big ] \, dx
\]
 attains its minimum at some  $\bar u\in W^{1,p}_{0}(\Omega)$ such that $I(\bar u)=0$, and 
 there exists  $\bar f \in L^q(\Omega; \bR^N)$ such that
\begin{eqnarray*}
\left\{
\begin{array}{lll}
\bar f (x) \in \bar \partial L(x, \nabla \bar u(x))&  \text{ a.e. } x \in \Omega\\
-{\rm div}(\bar f) = u^*.
\end{array}
\right.
\end{eqnarray*}
\end{pro}
\textbf{Proof.} Take $f_0 \in L^q(\Omega; \bR^N)$ with $-{\rm div} \big( f_0(x) \big) = u^*(x).$ Since  $L$ is an $\Omega$-dependent selfdual Lagrangian,  $M(x,a,b):=L(x,a,b+f_0(x))-\langle a, f_0(x) \rangle $ is also an $\Omega$-dependent selfdual Lagrangian on $\Omega \times \mathbb{R}^{N} \times \mathbb{R}^{N}.$  It follows from the above proposition that
\begin{eqnarray*}
F(v,v^*) := \inf_{ \substack{ f \in L^q(\Omega; \bR^N) \\ -{\rm div}(f) = v^* } } \S M \big( x,\nabla v(x),f(x) \big) \, dx
\end{eqnarray*}
is a selfdual Lagrangian on $W^{1,p}_{0}(\Omega) \times W^{-1,q}(\Omega).$ In view of the coercivity condition, Theorem \ref{one} applies and there exists $\bar u \in W^{1,p}_{0}(\Omega)$ such that
\[F(\bar u, 0)=\inf_{ \substack{ f \in L^q(\Omega; \bR^N) \\ -{\rm div}(f) = 0 } } \S M \big( x,\nabla \bar u(x),f(x) \big) \, dx=0. \]
Using again the coercivity condition, we get that the above infimum is attained at some $ f_1 \in L^q(\Omega; \bR^N)$  with ${\rm {\rm div}}( f_1) = 0.$ Thus,
\begin{eqnarray*}
0=F(\bar u, 0)&=& \S M \big( x,\nabla\bar u(x), f_1(x) \big) \, dx\\
&=& \S \big [L(x, \nabla \bar u(x),f_1(x)+f_0(x))-\langle \nabla \bar u(x), f_0(x)\rangle_{\bR^N} \big ] \, dx\\
&=& \S \big [L(x, \nabla \bar u(x),f_1(x)+f_0(x))-\langle \nabla \bar u(x), f_1(x)+f_0(x)\rangle_{\bR^N} \big ] \, dx.
\end{eqnarray*}
Taking into consideration that $L(x, \nabla \bar u(x),f_1(x)+f_0(x))-\langle \nabla \bar u(x), f_1(x)+f_0(x)\rangle_{\bR^N} \geq 0,$ we obtain that the latter is indeed zero, i.e.,
\[L(x, \nabla \bar u(x),f_1(x)+f_0(x))-\langle \nabla \bar u(x), f_1(x)+f_0(x)\rangle_{\bR^N}=0 \qquad \text{ for a.e. } x \in \Omega.\]
Setting $\bar f:= f_1+f_0$, we finally get that $\bar f(x) \in \bar \partial L (x, \nabla \bar u(x))$ for a.e. $x \in \Omega$ and that $-{\rm div}(\bar f) = u^*.$
\hfill $\square$

\subsection{Variational formula for the homogenized maximal monotone vector field}

Given a maximal monotone family ${\ss}$ in  $M_{\Omega, p}(\bR^N)$ that is $Q$-periodic for an  open non-degenerate parallelogram $Q$ in  $\bR^n$, its homogenization ${\ss}_{hom}$  is normally given by the non-variational formula (\ref{beta.hom}).
In this section, we shall give a variational formulation 
for the vector field ${\ss}_{hom}$  in terms of a suitably homogenized selfdual Lagrangian $L_{hom}$ derived from the $\Omega$-dependent selfdual Lagrangian associated to ${\ss}$. 

\begin{thm} \label{hom} Assume  ${\ss} \in M_{\Omega, p}(\bR^N)$ is $Q$-periodic and let $L$ be an   $\Omega$-dependent selfdual Lagrangian  such that ${\ss}(x, .)=\bar \partial L(x,.)$ given by Proposition \ref{rep}. If  the operator 
${\ss}_{hom}$ is given by (\ref{beta.hom}), then   ${\ss}_{hom}=\bar \partial L_{hom}$ where $L_{hom}$ is the selfdual Lagrangian on $\bR^N\times \bR^N$ given by 
\begin{equation}
L_{hom}(\xi,\eta)=\min_{\substack{\phi \in W^{1,p}_\# (Q) \\ g \in L_\#^q(Q; \bR^N)}}\frac{1}{|Q|}\int_Q L\big( x,\xi + \nabla \phi(x),\eta + g(x) \big) \dx.
\end{equation}
\end{thm}
The proof will follow from the following propositions. First, we 
show that the homogenized Lagrangian  $L_{hom }$
inherits many of the properties of the original $\Omega$-dependent Lagrangian $L$ such as convexity, boundedness and coercivity. 

\begin{pro}\label{l-infty} Assume $L$ is an $\Omega$-dependent Lagrangian on $\Omega \times \bR^N\times \bR^N$ satisfying (\ref{est200}) for some $p, q >1$. 
Then $L_{hom}$ is convex and lower semi continuous, and  for every $a^*,b^* \in \bR^n$,  
\begin{equation}\label{sd-hom}
L^*_{hom} (a^*,b^*)=\inf_{\substack{ \phi  \in W_{\#}^{1,q'}(Q)  \\ g \in L_{\#}^{p'}(Q; \bR^N) }}  \frac{1}{|Q|}\int_Q L^*\big( x,a^* + g(x), b^* +\nabla \phi (x) \big) \dx ,
\end{equation}
where $\frac{1}{p}+\frac{1}{p'}=1$ and $\frac{1}{q}+\frac{1}{q'}=1.$ Furthermore,
\begin{equation}\label{bd-hom}
C_0( |a|^p + |b|^q-1)\le L_{hom}(a,b) \le C_1( 1 + |a|^p + |b|^q) \qquad \text{ for all } a,b \in \bR^n.
\end{equation}
\end{pro}
The following gives the relation between the subdifferentials of $L_{hom}$ and of $L.$
\begin{pro}\label{SubDifAve}
 For each $a,b \in \bR^n$, the subdifferential map $\partial L_{hom}(a,b)$ is given by
\[
\partial L_{hom}(a,b) = \frac{1}{|Q|}\int_Q \partial L\big( y, a + \nabla \tilde{\phi}(y), b + \tilde{g}(y) \big) \dy,
\]
where $\tilde \phi  \in W_{\#}^{1,p}(Q)$ and   $\tilde g \in L_{\#}^{q}(Q; \bR^N)$ are such that 
\[
L_{hom}(a,b) = \frac{1}{|Q|}\int_Q  L\big( y, a + \nabla \tilde{\phi}(y), b + \tilde{g}(y) \big) \dy.
\]
\end{pro}
We need a few  preliminary facts. For each $1<r <\infty$, set
\[
E_r :=  \{ f= \nabla u \in  L^r(Q; \bR^N); \, \mbox{for some }u \in W^{1,r}_\#(Q)\}
\]
and
\[
E_r + \bR^n:= \{ f + \eta \; : \; f \in E_r, \, \eta \in \bR^n \}.
\]
The Poincar\'e-Wirtenger inequality which states that for $D$ bounded open and convex, there exists  $K:=K(r,D)>0$ such that 
\[
\hbox{$\|u- \frac{1}{|D|}\int_D u\|_{_{L^r(D)}}\le K \|\nabla u\|_{W^{1,r}(D)}$ for every $u \in W^{1,r}(D), $}
\] 
 implies that  $E_r + \bR^n$ is a convex weakly closed subset of $L^r(Q; \bR^N).$
The  indicator function of $E_r + \bR^n$,
\[
\chi_{E_r + \bR^n}(f) = \left\{ \begin{array}{lc} 0 & f \in E_r + \bR^n, \\ +\infty & \qquad \qquad \quad \quad  f \in L^r(Q; \bR^N)\setminus (E_r + \bR^n),  \end{array}\right.
\]
is therefore convex and lower semi-continuous in $L^r(Q; \bR^N)$. Assuming that $r'$ is the conjugate of $r$, i.e., $\frac{1}{r}+\frac{1}{r'}=1,$ define
\[E_{r'}^\perp := \big \{ g \in L^{r'}(Q; \bR^N); \, \S \langle f(x),g(x) \rangle_{\bR^N} \, dx= 0 \quad \text{ for all } f \in  E_r + \bR^n \big \}.\]
The Fenchel-Legendre dual  $\chi^*_{E_r + \bR^n}$ of  $\chi_{E_r + \bR^n}$ is then given by,
\begin{eqnarray*}
\chi^*_{E_r + \bR^n} (g) &=& \sup_{f \in L^r(Q; \bR^N)} \Big\{\int_Q \langle f(x), g(x) \rangle_{\bR^N} \, dx-\chi_{E + \bR^n}(f)\Big \}\\
&=& \sup_{f \in E_r + \bR^n} \int_Q \langle f(x), g(x) \rangle_{\bR^N} \, dx=\chi_{E_{r'}^\perp} (g),
\end{eqnarray*}
 for all $g \in L^{r'}(Q; \bR^N)$. Also due to the convexity and lower semi-continuity of $\chi_{E_r + \bR^n}$ one has $\chi^*_{E_{r'}^\perp}=\chi_{E_r + \bR^n}.$ Similarly one can deduce  that,
\[
\chi^*_{E_{r'}^\perp + \bR^n} (f) = \chi_{E_r} (f)
\]
 for all $ f \in L^r(Q; \bR^N).$ Note also that $E_r$ is the isometric image of $W_{\#}^{1,r}(Q)$ by $\nabla$ and $E_{r}^\perp=L_{\#}^r(Q; \bR^N).$\\

\textbf{Proof of Proposition \ref{l-infty}.} We first prove (\ref{sd-hom}).
Fix $(a^*,b^*) \in \bR^n \times \bR^n$ and write
\begin{eqnarray*}
L^*_{hom} (a^*,b^*) &=& \sup_{(a,b) \in \bR^n \times \bR^n} \{ \langle a,a^* \rangle_{\bR^N} + \langle b,b^* \rangle_{\bR^N} - L_{hom} (a,b)\} \\
&=& \sup_{\substack{(a,b) \in \bR^n \times \bR^n \\ (\phi,g) \in W_{\#}^{1,p}(Q) \times L_{\#}^q(Q; \bR^N) }} \frac{1}{|Q|} \int_Q \Big[ \langle a,a^* \rangle_{\bR^N} + \langle b,b^* \rangle_{\bR^N} - L\big( x,a + \nabla \phi(x), b + g(x) \big) \Big] \dx \\
&=& \sup_{\substack{(a,b) \in \bR^n \times \bR^n \\ (f,g) \in E_p \times E_q^\perp }}  \frac{1}{|Q|} \int_Q \Big[ \langle a,a^* \rangle_{\bR^N} + \langle b,b^* \rangle_{\bR^N} - L\big( x,a + f(x), b + g(x) \big) \Big] \dx .
\end{eqnarray*}

Setting $A(x) = a + f(x)$, $B(x) = b + g(x)$ and substituting above we have
\begin{eqnarray*}
L^*_{hom} (a^*,b^*) &=& \sup_{\substack{A \in E_p + \bR^n \\ B \in E_q^\perp + \bR^n }}  \frac{1}{|Q|} \int_Q \Big[ \langle A,a^* \rangle_{\bR^N} + \langle B,b^* \rangle_{\bR^N} - L\big( x,A(x), B(x) \big) \Big] \dx \\
  &=& \sup_{\substack{A \in L^p(\Omega; \bR^n ) \\ B \in L^q(\Omega; \bR^n) }} \Big \{ \frac{1}{|Q|} \int_Q \Big[ \langle A,a^* \rangle_{\bR^N} + \langle B,b^* \rangle_{\bR^N} - L\big( x,A(x), B(x) \big) \Big] \dx \\&&-\chi_{E_p + \bR^n}(A) - \chi_{E_q^\perp + \bR^n}(B) \Big\}.
\end{eqnarray*}
Now using the fact that the Fenchel dual  of a sum is their  \emph{inf-convolution}, we obtain
\begin{eqnarray*}
L^*_{hom} (a^*,b^*) &=&
\inf_{\substack{ f \in L^{q'}(Q; \bR^n) \\ g \in L^{p'}(Q; \bR^n) }} \big \{ \frac{1}{|Q|}\int_Q L^*\big( x,a^* - g(x), b^* - f(x) \big) \dx + \chi_{E_{p'}^\perp}(g) + \chi_{E_{q'}} (f) \big \}\\
&=& \inf_{\substack{ f \in E_{q'} \\ g \in E_{p'}^\perp }}  \frac{1}{|Q|}\int_Q L^*\big( x,a^* - g(x), b^* - f(x) \big) \dx \\
&=& \inf_{\substack{ \phi  \in W_{\#}^{1,q'}(Q)  \\ g \in L_{\#}^{p'}(Q; \bR^N) }}\frac{1}{|Q|}\int_Q L^*\big( x,a^* + g(x), b^* +\nabla \phi (x) \big) \dx.
\end{eqnarray*}
This proves (\ref{sd-hom}), which then implies that 
$L^{**}_{hom}=L_{hom}$ and therefore $L_{hom}$ is convex and lower semi-continuous.\\
We now prove estimate  (\ref{bd-hom}). In fact,  the upper bound  simply follows from 
\[
L_{hom} (a,b) \le \frac{1}{|Q|} \int_Q L(x,a,b) \dx \le C_1(|a|^p + |b|^q + 1).
\]
 For the lower bound, note first that since  $C_0(|a|^p + |b|^q - 1) \le L(x,a,b)$ for all $a,b \in \bR^N,$ it follows that
\[ L^*(x,a,b)  \le \frac{C_0(p-1)}{(C_0p)^{p'}}|a|^{p'} + \frac{C_0(q-1)}{(C_0q)^{q'}}|b|^{q'} +C_0  \quad \text{ for all } \quad a,b \in \bR^N.\]
On then get from  (\ref{sd-hom}) that
\[
L^*_{hom} (a,b) \le \frac{1}{|Q|} \int_Q L^*(x,a,b) \dx \le \frac{C_0(p-1)}{(C_0p)^{p'}}|a|^{p'} + \frac{C_0(q-1)}{(C_0q)^{q'}}|b|^{q'} +C_0  \quad \text{ for all } \quad a,b \in \bR^N,
\]
from which we get that  $
L_{hom} (a,b)= L^{**}_{hom} (a,b)\ge C_0(|a|^p + |b|^q -1)$ 
for all $ a,b \in \bR^N.$
\hfill $\square$\\

\textbf{Proof of Proposition \ref{SubDifAve}.} Setting $
A(a,b) := \frac{1}{|Q|}\int_Q \partial L\big( y, a + \nabla \tilde{\phi}(y), b + \tilde{g}(y) \big) \dy, 
$ we shall first show that $A \subset \partial L_{hom}$. For that  consider $(a_1,b_1) \in \bR^N \times \bR^N, \phi  \in W_{\#}^{1,p}(Q)$ and $  g \in L_{\#}^{q}(Q; \bR^N).$ From the convexity of $L$:
\begin{eqnarray*}
L\big( y, a_1 + \nabla {\phi}(y), b_1 + {g}(y) \big)&\ge& L\big( y, a + \nabla \tilde{\phi}(y), b + \tilde{g}(y) \big)\\
&&+\langle \partial_1 L\big( y, a + \nabla\tilde{\phi}(y), b + \tilde{g}(y) \big) ,  a_1+ \nabla {\phi}(y) - a - \nabla\tilde{\phi}(y) \rangle_{\bR^N}\\
&&+\langle \partial_2 L\big( y, a + \nabla\tilde{\phi}(y), b + \tilde{g}(y) \big) ,  b_1 + {g}(y) - b- \tilde{g}(y)  \rangle_{\bR^N}.
\end{eqnarray*}
 Averaging the above on $Q$ implies that
\[
\frac{1}{|Q|}\int_Q L\big( y, a_1 + \nabla {\phi}(y), b_1 + {g}(y) \big) \dy \ge L_{hom}( a, b ) + \langle  A ( a , b ) , ( a_1 - a , b_1 - b) \rangle_{\bR^N \times \bR^N}, 
\]
from which we get
\[
L_{hom}( a_1, b_1 )  \ge L_{hom}( a, b ) + \langle  A ( a , b ) , ( a_1 - a , b_1 - b ) \rangle_{\bR^N \times \bR^N}.
\]
This implies that $A \subset \partial L_{hom}.$  To prove the reverse inclusion, let $(d,c)$ be in $\partial L_{hom} (a,b).$ Since  $L_{hom}$ is convex and lower semi-continuous,  we have
\[L_{hom}(a,b)+L^*_{hom}(d,c)= \langle a,d\rangle_{\bR^N }+\langle b,c\rangle_{\bR^N }.\]
It follows from Proposition \ref{l-infty} that there exist $\phi  \in W_{\#}^{1,q'}(Q) $ and $ g \in L_{\#}^{p'}(Q; \bR^N)$ such that
\begin{equation*}
L^*_{hom} (a^*,b^*)= \frac{1}{|Q|}\int_Q L^*\big( x,a^* + g(x), b^* +\nabla \phi (x) \big) \dx, 
\end{equation*}
and therefore
\[
\frac{1}{|Q|}\int_Q  L\big( y, a + \nabla \tilde{\phi}(y), b + \tilde{g}(y) \big) \dy+\frac{1}{|Q|}\int_Q L^*\big( x,d + g(x), c +\nabla \phi (x) \big) \dx=\langle a,d\rangle_{\bR^N }+\langle b,c\rangle_{\bR^N }.
\]
On the other hand,
 \[\langle a,d\rangle_{\bR^N }+\langle b,c\rangle_{\bR^N }=\frac{1}{|Q|}\int_Q \langle a+\nabla \tilde{\phi}(y),d+g(y)\rangle_{\bR^N }+\frac{1}{|Q|}\int_Q \langle b+\tilde{g}(y),c+\nabla \phi (y)\rangle_{\bR^N } \, dy,\]
which together with the previous equality yield
\[\int_Q \big [ L\big( y, a + \nabla \tilde{\phi}(y), b + \tilde{g}(y) \big) +L^*\big( y,d + g(y), c +\nabla \phi (y) \big)-\langle a+\nabla \tilde{\phi}(y),d+g(y)\rangle_{\bR^N }- \langle b+\tilde{g}(y),c+\nabla \phi (y)\rangle _{\bR^N } ]\, dy=0.\]
Taking into account that the integrand is non-negative we obtain
\[ L\big( y, a + \nabla\tilde{\phi}(y), b + \tilde{g}(y) \big) +L^*\big( y,d + g(y), c +\nabla \phi (y) \big)-\langle a+\nabla \tilde{\phi}(y),d+g(y)\rangle_{\bR^N }- \langle b+\tilde{g}(y),c+\nabla \phi (y)\rangle_{\bR^N }=0\]
for almost all $y \in Q.$ This implies that
\[(d + g(y), c +\nabla \phi (y))\in \partial L\big( y, a + \nabla \tilde{\phi}(y), b + \tilde{g}(y) \big) \qquad \text{ a.e. } y \in Q. \]
Integrating the above over $Q$ implies that
\[(d , c )\in \frac{1}{|Q|}\int_Q\partial L\big( y, a + \nabla \tilde{\phi}(y), b + \tilde{g}(y) \big),
 \]
which completes the proof.
\hfill $\square$\\

\textbf{Proof of Theorem \ref{hom}}  Let $\eta  \in \bar \partial L_{hom}(\xi)$ in such a way that 
$
L_{hom}(\xi,\eta)=\langle \xi, \eta \rangle_{\bR^N}.
$ 
From the definition of $L_{hom}$, we have
\[L_{hom}(\xi,\eta)=\min_{\substack{\phi \in W^{1,p}_\# (Q) \\ g \in L_\#^q(Q; \bR^N)}}\frac{1}{|Q|}\int_Q L\big( x,\xi + \nabla \phi(x),\eta + g(x) \big) \dx.\]
From the coercivity assumptions on $L$, it follows that there exist $\phi \in W^{1,p}_\# (Q)$ and $g \in L_\#^q(Q; \bR^N)$ such that
\[L_{hom}(\xi,\eta)=\frac{1}{|Q|}\int_Q L\big( x,\xi + D\phi(x),\eta + g(x) \big) \dx.\]
Hence
\begin{eqnarray*}
0&=&L_{hom}(\xi,\eta)-\langle \xi, \eta \rangle_{\bR^N}\\
&=&\frac{1}{|Q|}\int_Q L\big( x,\xi + \nabla \phi(x),\eta  + g(x) \big) \dx-\langle \xi, \eta \rangle_{\bR^N}\\
&=&\frac{1}{|Q|}\int_Q\big [ L\big( x,\xi + \nabla \phi(x),\eta  + g(x) \big) -\langle \xi + \nabla \phi(x), \eta + g(x)\rangle_{\bR^N} \big ]\dx,
\end{eqnarray*}
and since the integrand in non-negative we obtain
\[L\big( x,\xi + \nabla \phi(x),\eta  + g(x) \big) -\langle \xi + \nabla \phi(x), \eta + g(x)\rangle_{\bR^N}=0 \,\, \text{ for a.e. } x \in Q,
\]
from which we have
\[\eta + g(x) \in \bar \partial L(x,\xi + \nabla \phi(x))={\ss}(x, \xi  + \nabla \phi(x)) \]
and finally
$\eta = \int_Q(\eta +g(x))\, dx.$ This implies that $\bar \partial L_{hom}\subset {\ss}_{hom} $ and the equality follows since $\bar \partial L_{hom} $ is itself a maximal monotone operator. \hfill $\square$

\section{A variational approach to homogenization}

We start by studying the homogenization of a  class of Lagrangians that is more general  than the one introduced in Proposition \ref{lift}. We shall then apply this result to deduce Theorem \ref{main} announced in the introduction.

\subsection{The homogenization  of general Lagrangians on $W^{1,p}(\Omega) \times L^q(\Omega; \bR^N)$}\label{HomogResults}

The following homogenization result does not require  the $\Omega$-dependent Lagrangian $L$ to be selfdual nor that the exponents $p$ and $q$ to be conjugate. 

\begin{thm}\label{gam-lim} Let $\Omega$ be a regular bounded domain and $Q$   an open non-degenerate parallelogram in  $\bR^n$. Let  $L: \Omega \times \mathbb{R}^{N} \times \mathbb{R}^{N}\to \R $ be  an  $\Omega$-dependent Lagrangian such that: \\
(1) For each $a,b \in \bR^N $ the function
$x \to L(x,a,b)$  is $Q$-periodic.\\
(2)  There exist constants $C_0,C_1 \ge 0$ and exponents $p,q>1$ such that for every $x\in \Omega$,
\[
C_0(|a|^p + |b|^q - 1) \le L(x,a,b) \le C_1(|a|^p + |b|^q + 1).
\]
 Let $\{G_\e; \e >0\}$ be the family of functionals on  $W^{1,p}(\Omega) \times L^q(\Omega;\bR^N)$ defined by
\[
G_\e (u,\tau) := \inf_{\substack{f \in L^q(\Omega;\bR^N) \\ {\rm div} f = 0}} \S L\big(\frac{x}{\e},\nabla u(x),\tau(x) + f(x)\big) \dx,
\]
and set 
\begin{equation}\label{hom-int}
L_{hom}(a,b) := \min_{\substack{\phi \in W^{1,p}_\# (Q) \\ g \in L_\#^q(Q; \bR^N)}}\frac{1}{|Q|}\int_Q L\big( x,a + \nabla \phi(x),b + g(x) \big) \dx.
\end{equation}
Equip $L^q(\Omega;\bR^N)$ with the following topology denoted by ${\mathcal T},$
\[\tau_n\to \tau \text{ for } {\mathcal T} \quad \text{\rm if and only if } \quad \tau_n\to \tau \text{ weakly in }  L^q(\Omega;\bR^N) \text{ and } {\rm div} (\tau_n)\to {\rm div} (\tau) \text{ strongly in } W^{-1,q}(\Omega),
\]
There exists then a Lagrangian $G_{hom}$ 
on $W^{1,p}(\Omega) \times L^q(\Omega;\bR^N)$ that is a $\Gamma$-limit of $\{G_\e; \e >0\}$ as ${\e \to 0}$. Moreover, $G_{hom}$ is given by the formula
\begin{equation}\label{hom-lag}
G_{hom} (u,\tau) := \inf_{\substack{f \in L^q(\Omega;\bR^N)\\ {\rm div} f = 0}} \S L_{hom}\big( \nabla u(x),\tau (x) + f(x)\big) \dx,
\end{equation}
\end{thm}

\begin{remark} \rm 
Note that when the Lagrangian $L$ is independent of the third  variable, i.e.,
\[L(x,a,b)=\phi(x,a) \quad \text{ for all } (x,a,b) \in \Omega \times \mathbb{R}^{N} \times \mathbb{R}^{N},\]
for some function $\phi:\Omega \times \mathbb{R}^{N} \to \mathbb{R},$ this homogenization  problem is completely understood. Also,  when the Lagrangian $L$ is independent of the second   variable then this problem can be dealt using the bi-continuity of the Fenchel dual (see for instance \cite{Att-book,D-M-S-Mosco}). The proof for the general Lagrangians consists of two parallel  parts corresponding to each of these variables and  should be done simultaneously for both. The part regarding the first variable is rather  standard and the same argument can be found for instance in \cite{Att-book}.
\end{remark}
The proof of Theorem \ref{gam-lim} will follow from the following two lemmas.
\begin{lem}\label{proof1} For  any $(u,\tau) \in W^{1,p}(\Omega) \times L^q(\Omega;\bR^N)$, there exists a  sequence $(u_\e,\tau_\e)\in W^{1,p}(\Omega) \times L^q(\Omega;\bR^N)$ such that $u_\e\to u$ strongly in $L^p(\Omega)$, $\tau_\e \to \tau$ strongly in  $L^q(\Omega;\bR^N)$ and
\[
\limsup_{\e \to 0} G_\e (u_\e,\tau_\e) \le G_{hom}(u,\tau ).
\]
\end{lem}
\begin{lem}\label{proof2} Let $f \in L^q(\Omega;\bR^N)$  with ${\rm div}(f)=0.$ For  any $(u,\tau) \in W^{1,p}(\Omega) \times L^q(\Omega;\bR^N)$ and any   sequence $(u_\e,\tau_\e)$ such that $u_\e\to u$ strongly in $L^p(\Omega)$ and  $\tau_\e \to \tau$ with the ${\mathcal T} $-topology  in  $L^q(\Omega;\bR^N)$,  we have
\[
\liminf_{\e \to 0} \S L(\frac{x}{\e},\nabla u_\e(x),\tau_\e(x) + f(x) ) \dx \ge \S L_{hom}\big( \nabla u(x),\tau(x) + f(x) \big) \dx.
\]
\end{lem}
We first show how Theorem \ref{gam-lim}} follows from the two lemmas above. 

The limsup property in the definition of $\Gamma$-convergence readily follows from Lemma \ref{proof1}.  For the liminf property we must show that for any $(u,\tau) \in W^{1,p}(\Omega) \times L^q(\Omega; \bR^N)$ and any sequence $\{ (u_\e,\tau_\e) \} \subset W^{1,p}(\Omega) \times L^q(\Omega; \bR^N)$ such that
\[
u_\e \to u  \text{ strongly in } L^p(\Omega) \quad \text{ and } \quad \tau_\e \to \tau \quad \text{  in the}\,\,   {\mathcal T}-\text{topology}, 
\]
we have that
\[
\liminf_{\e \to 0} G_\e (u_\e,\tau_\e) \ge G_{hom} (u,\tau ).
\]
By  Lemma \ref{proof2} we have
 \[
\liminf_{\e \to 0} \S L(\frac{x}{\e}, \nabla u_\e,\tau_\e + f ) \dx \ge \S L_{hom}\big( \nabla u,\tau + f \big) \dx,
\]
for every $f \in L^q(\Omega;\bR^N)$  with ${\rm div}(f)=0.$ Since

\[
\inf_{\substack{f \in L^q(\Omega; \bR^N) \\ {\rm div} f = 0}} \liminf_{\e \to 0} \S L(\frac{x}{\e},\nabla u_\e,\tau_\e + f ) \dx = \liminf_{\e \to 0} \inf_{\substack{f \in L^q(\Omega; \bR^N) \\ {\rm div} f = 0}} \S L(\frac{x}{\e},\nabla u_\e,\tau_\e + f ) \dx,
\]
we obtain that $
\liminf_{\e \to 0} G_\e (u_\e,\tau_\e) \ge G_{hom} (u,\tau ),$ 
as desired. \hfill $\square$\\

\textbf{Proof of Lemma  \ref{proof1}.} Note that without  loss of generality we may assume $L\geq 0.$
 Assume first that $u$ is  an affine function and $\tau$ is constant on $\Omega$, that is
 \[
 \hbox{$
u(x) = \langle a, x \rangle + \alpha$
and $
\tau (x) = b$, }
\]
for some $a$ and $b$ in $\bR^n$ and $\alpha \in \bR$.  Fix $\eta \in \bR^n$ and let  $\tilde{\phi}$ and $\tilde{g}$ to be the minimizers on the formula for $L_{hom}$ given by (\ref{hom-int}):
\begin{equation}\label{min200}
L_{hom} (a,b+\eta) = \frac{1}{|Q|} \int_Q L\big( x, a + \nabla \tilde{\phi}(x), b + \eta + \tilde{g}(x) \big).
\end{equation}
Define 
\[
\hbox{$u_\e(x) := u(x) + \e \tilde{\phi}(\frac{x}{\e})$\quad 
and \quad $
\tau_\e(x) := \tau.$}
\]
Note that by Lemma \ref{ext2} in the Appendix, $\tilde{g}$ can be extended by periodicity to an element of $L_{loc}^{q}(\bR^N; \bR^N),$ still denoted by $\tilde{g}$ such that ${\rm div} (\tilde{g})=0.$ It follows that
\begin{eqnarray*}
\limsup_\e G_\e (u_\e,\tau_\e) & = &
\limsup_\e \inf_{ \substack{f \in L^q(\Omega; \bR^N) \\ {\rm div} f = 0} } \S L\big(\frac{x}{\e},a + \nabla \tilde{\phi}(\frac{x}{\e}), b + f(x) \big) \dx\\
 & \le & \inf_{ \substack{f \in L^q(\Omega; \bR^N) \\ {\rm div} f = 0} } \limsup_\e \S L\big(\frac{x}{\e},a + \nabla \tilde{\phi}(\frac{x}{\e}), b + f(x) \big) \dx\\
 & \le & \limsup_\e \S L\big(\frac{x}{\e},a + \nabla \tilde{\phi}(\frac{x}{\e}), b + \eta + \tilde{g}(\frac{x}{\e}) \big) \dx.
\end{eqnarray*}
By Lemma \ref{Riem-Leb} of the Appendix we have as $\e \to 0$,
\[
\S L\big(\frac{x}{\e},a + \nabla \tilde{\phi}(\frac{x}{\e}), b + \eta + \tilde{g}(\frac{x}{\e}) \big) \dx \to
\frac{|\Omega|}{|Q|}\int_Q L\big(y,a + \nabla \tilde{\phi}(y), b + \eta + \tilde{g}(y) \big) \dy.
\]
It then follows from (\ref{min200}) that
\[
\limsup_{\e \to 0} G_\e (u_\e,\tau_\e) \le |\Omega| L_{hom} (a,b+\eta),
\]
and since $\eta$ is arbitrary, we have that
\[
\limsup_{\e \to 0}  G_\e (u_\e,\tau_\e) \le \inf_{\eta \in \bR^n}|\Omega| L_{hom} (a,b+\eta).
\]
By Lemma \ref{II} of the Appendix we have
\begin{eqnarray*}
\inf_{\substack{f \in L^q(\Omega; \bR^N) \\ {\rm div} f = 0}} \S L_{hom}\big(a,b + f(x)\big) \dx  \geq  \inf_{\eta \in \bR^n}|\Omega| L_{hom} (a,b+\eta),
\end{eqnarray*}
and thus we conclude, as desired
\[
\limsup_{\e \to 0}  G_\e (u_\e,\tau_\e) \le \inf_{\substack{f \in L^q(\Omega; \bR^N) \\ {\rm div} f = 0}} \S L_{hom}\big(a,b + f(x)\big) \dx = G_{hom}(u,\tau).
\]
Assume now that $u$ is a piecewise affine function and $\tau$ is a piecewise constant function  on $\Omega$, that is  for $\{\hat{\Omega}_j\}_{j \in I_1}$ and $\{\tilde{\Omega}_k\}_{k \in I_2}$, both finite polyhedral partitions of $\Omega$, we have
\[
\hbox{$u(x) = \langle a_j,x \rangle + \alpha_j$ for $ x \in \hat{\Omega}_j$ \quad 
and \quad
 $\tau(x) =  b_k$ for $x \in \tilde{\Omega}_k,$}
\]
for fixed $a_j \in \bR^n$ and $b_k \in \bR^n$ and constants $\alpha_j.$
 By considering non-empty intersections  $\hat{\Omega}_j\cap\tilde{\Omega}_k$  and re-indexing them, we can consider $\{\Omega_i\}_{i \in I}$ a polyhedral partition of $\Omega$ such that
\[
\hbox{$ u(x) = \langle a_i,x \rangle + \alpha_i$ for $ x \in \Omega_i$\quad 
and \quad $\tau(x) =  b_i$ for $ x \in \Omega_i$.}
\]
Analogous to what was done previously, fix  $\{\eta_i\} \subset \bR^N$ and let $\tilde{\phi}_i$ and $\tilde{g}_i$ be  such that
\[
L_{hom} (a_i,b_i+\eta_i) = \frac{1}{|Q|} \int_Q L\big( x, a_i + \nabla \tilde{\phi}_i(x), b_i + \eta_i + \tilde{g}_i(x) \big)\dx, 
\]
and set $u^i_\e(x) := u(x) + \e \tilde{\phi}_i(\frac{x}{\e})$.

Unfortunately, we cannot consider $u_\e$ as the piecewise construction of the above functions, as the $\phi_i$ won't necessarily match along the interface between the $\Omega_i$ and thus will not in general be a function in $W^{1,p}(\Omega)$. This can be remedied by the following standard construction (see for instance \cite{Att-book}):  
Let $\Sigma$ be the interface set between the $\Omega_i$, and define for $\delta > 0$, $
\Sigma_\delta := \{ x\in\Omega \; : \; d(x,\Sigma) \le \delta \}.$
Consider a smooth function $\Psi_\delta$ so that
\[
\Psi_\delta (x)=\left\{
\begin{array}{cl}
1 & x\in\Sigma_\delta \\
0 & x \in \Omega \setminus \Sigma_{2\delta}, 
\end{array}
\right.
\]
and define
\[
\hbox{$u^\delta_\e(x) := \big( 1 - \Psi_\delta(x) \big) u^i_\e (x) + \Psi_\delta (x)u(x) \;\mbox{ for } x\in\Omega_i$\quad 
and \quad $
\tau_\e := \tau.$}
\]
 It can be checked that the function $u^\delta_\e$ lies in $W^{1,p}(\Omega)$.
 Note that by Lemma \ref{ext2} of the Appendix, each  $\tilde{g_i}$ can be extended by periodicity to an element of $L_{loc}^{q}(\bR^N; \bR^N),$ still denoted by $\tilde{g_i}$ such that ${\rm div} (\tilde{g_i})=0.$ Thus ${\rm div} (\eta_i+\tilde{g_i}(\frac{x}{\e}))=0$ on $\bR^N$ and in particular on $\Omega_i\setminus \Sigma_{\delta}. $ Define $f_{\e, \delta}(x)=\eta_i+\tilde{g_i}(\frac{x}{\e})$ on $\Omega_i\setminus \Sigma_{\delta}. $ One can also extend (using Theorem 2.5 and Corollary 2.8 in \cite{Gr}) $f_{\e, \delta}$ to an element in $L^q(\Omega; \bR^N),$ still denoted by $f_{\e, \delta}$ such that $\|f_{\e, \delta}\|_{L^q(\Omega; \bR^N)}$ is bounded and ${\rm div}(f_{\e, \delta})=0.$
 Take now any $0 < t < 1$, then 
 \begin{eqnarray*}
G_\e (tu^\delta_\e,\tau_\e) &=&
 \inf_{\substack{f \in L^q(\Omega; \bR^N) \\ {\rm div} f = 0}}  \int_{\Omega} L\Big(\nabla tu^\delta_\e,\tau_\e + f\Big) \dx \\ &\leq & \int_{\Omega} L\Big(\nabla tu^\delta_\e,\tau_\e + f_{\e, \delta}\Big) \dx\\ &=&
  \sum_i \int_{\Omega_i\setminus \Sigma_{\delta}} L\Big(\frac{x}{\e},t\big(1-\Psi_\delta\big)\nabla u^i_\e + t\Psi_\delta \nabla u + (1-t)\frac{t}{(1-t)}(u - u^i_\e)\nabla \Psi_\delta,b_i +\eta_i+\tilde{g_i}(\frac{x}{\e}) \Big) \dx\\&&+\int_{\Sigma_{\delta}} L\Big(\nabla tu^\delta_\e,\tau_\e + f_{\e, \delta}\Big) \dx
\end{eqnarray*}

Since $L$ is convex in the middle variable and since $t( 1-\Psi_\delta) + t\Psi_\delta +(1-t) = 1$, we obtain
\begin{eqnarray*}
G_\e (tu^\delta_\e,\tau_\e) &\le&
  \sum_i \int_{\Omega_i\setminus \Sigma_{\delta}} t( 1-\Psi_\delta) L\Big(\frac{x}{\e},a_i + \nabla \tilde{\phi}_i(\frac{x}{\e}),b_i+\eta_i+\tilde{g_i}(\frac{x}{\e}) \Big)\dx \\
&& + \sum_i \int_{\Omega_i\setminus \Sigma_{\delta}} (1-t)L\Big(\frac{x}{\e},\frac{t}{(1-t)}(u - u^i_\e)\nabla \Psi_\delta,b_i+\eta_i+\tilde{g_i}(\frac{x}{\e} \Big)\dx\\&& + \int_{\Sigma_{2\delta}\setminus \Sigma_{\delta}} t\Psi_\delta L\Big(\frac{x}{\e},\nabla u,b_i + \eta_i+\tilde{g_i}(\frac{x}{\e}\Big) \dx \\ && +\int_{\Sigma_{\delta}} L\Big(\nabla tu^\delta_\e,\tau_\e + f_{\e, \delta}\Big) \dx.
\end{eqnarray*}

For the first term on the right hand side of this inequality we have
\[
\int_{\Omega_i\setminus \Sigma_{\delta}} t( 1-\Psi_\delta)  L\Big(\frac{x}{\e},a_i + \nabla \tilde{\phi}_i(\frac{x}{\e}),b_i+\eta_i+\tilde{g_i}(\frac{x}{\e})  \Big)  \dx \le \int_{\Omega_i\setminus \Sigma_{\delta}} L\Big(\frac{x}{\e},a_i + \nabla \tilde{\phi}_i(\frac{x}{\e}),b_i+\eta_i+\tilde{g_i}(\frac{x}{\e})  \Big)\dx.
\]
Using the boundedness of $L$ we get the following estimate for the second term,
\begin{eqnarray*}
\int_{\Omega_i\setminus \Sigma_{\delta}} (1-t)L\Big(\frac{x}{\e},\frac{t}{(1-t)}(u - u^i_\e)\nabla \Psi_\delta,b_i+\eta_i+\tilde{g_i}(\frac{x}{\e})  \Big)\dx &\le & C_1 (1-t) \int_{\Omega_i\setminus \Sigma_{\delta}} \Big(|\frac{t}{(1-t)}(u - u^i_\e)\nabla \Psi_\delta|^p \\&&+ |b_i+\eta_i+\tilde{g_i}(\frac{x}{\e}) |^q +1\Big)\dx,
\end{eqnarray*}
and similarly
\begin{eqnarray*}
\int_{\Sigma_{2\delta}\setminus \Sigma_{\delta}} t\Psi_\delta L\Big(\frac{x}{\e},\nabla u,b_i + \eta_i+\tilde{g_i}(\frac{x}{\e})  \Big) \dx \le  C_1 \int_{\Sigma_{2\delta}\setminus \Sigma_{\delta}}\Big(1+|\nabla u|^p + |b_i + \eta_i+\tilde{g_i}(\frac{x}{\e}) |^q \Big) \dx,
\end{eqnarray*}
as well as
\[
\int_{\Sigma_{\delta}} L\Big(\nabla tu^\delta_\e,\tau_\e + f_{\e, \delta}\Big) \dx \le C_1 \int_{\Sigma_{\delta}}\Big(1+|\nabla tu^\delta_\e|^p + |\tau_\e + f_{\e, \delta} |^q \Big) \dx.
\]
It then follows that
\begin{eqnarray*}
G_\e (tu^\delta_\e,\tau_\e)
&\le& \sum_i
\int_{\Omega_i\setminus \Sigma_{\delta}} L\Big(\frac{x}{\e},a_i + \nabla \tilde{\phi}_i(\frac{x}{\e}),b_i + \eta_i+\tilde{g}_i(\frac{x}{\e}) \Big)\dx \\&&
+ C_1 (1-t) \sum_i\int_{\Omega_i\setminus \Sigma_{\delta}} \Big( |\frac{t}{(1-t)}(u - u^i_\e)\nabla \Psi_\delta|^p + |b_i+ \eta_i+\tilde{g}_i(\frac{x}{\e})|^q+1 \Big)\dx  \\&&
+ C_1 \int_{\Sigma_{2\delta}\setminus \Sigma_{\delta}}\Big(1+|\nabla u|^p + |b_i +  \eta_i+\tilde{g}_i(\frac{x}{\e})|^q \Big) \dx\\ &&+
C_1 \int_{\Sigma_{\delta}}\Big(1+|\nabla tu^\delta_\e|^p + |\tau_\e + f_{\e, \delta} |^q \Big) \dx.
\end{eqnarray*}
By  taking   $\limsup_{\e \to 0} $ on both sides and considering  $u^i_\e \to u$ on $L^p(\Omega_i)$, and then letting $t \to 1$ and $\delta \to 0$ we finally get,
\begin{equation}\label{lim-lim}
\limsup_{\substack{t \to 1 \\ \delta \to 0}}\limsup_{\e \to 0}  G_\e (tu^\delta_\e,\tau_\e) \le \sum_i \frac{|\Omega_i|}{|Q|}\int_{Q} L\big( x, a_i + \nabla \tilde{\phi}_i(x), b_i + \eta_i + \tilde{g}_i(x) \big)\dx.
\end{equation}
Also note that,
\[
\sum_i \frac{|\Omega_i|}{|Q|}\int_Q L\big( x, a_i + \nabla \tilde{\phi}_i(x), b_i + \eta_i + \tilde{g}_i(x) \big)\dx = \sum_i |\Omega_i| L_{hom}(a_i,b_i+\eta_i).
\]
 A diagonalization argument yields from limit \eqref{lim-lim} the existence of some $t(\e)$ and $\delta(\e)$ such that $t(\e) \to 1$ and $\delta(\e) \to 0$ as $\e \to 0.$ Defining
$
u_\e : = t(\e)u^{\delta(\e)}_\e,
$
we obtain 
\[
\limsup_{\e \to 0}  G_\e (u_\e,\tau_\e) \le \sum_i |\Omega_i| L_{hom}(a_i,b_i+\eta_i),
\]
and since the $\{\eta_i\}$ is  arbitrary one has
\[
\limsup_{\e \to 0}  G_\e (u_\e,\tau_\e) \le \sum_i |\Omega_i| \inf_{\eta_i \in \bR^n }
L_{hom}(a_i,b_i+\eta_i).
\]
 Now we use Lemma \ref{II} of the Appendix
to   obtain
\[
\sum_i |\Omega_i| \inf_{\eta_i \in \bR^n }
L_{hom}(a_i,b_i+\eta_i) \le \inf_{\substack{f \in L^q(\Omega; \bR^N) \\ {\rm div} f = 0}} \S L_{hom}\big(\nabla u(x),\tau(x) + f(x)\big) \dx,
\]
from which we get
$
\limsup_{\e \to 0}  G_\e (u_\e,\tau_\e) \le G_{hom}(u,\tau).
$

 Finally, consider any $(u,\tau) \in W^{1,p}(\Omega) \times L^q(\Omega; \bR^N)$. There exists then a sequence $\{ u_n \}$ of piecewise affine functions and a sequence $\{ \tau_n \}$ of piecewise constant functions such that $(u_n,\tau_n) \to (u,\tau)$. By Proposition \ref{l-infty}, the function $G_{hom}$ are  continuous, so we also have
\[
\lim_n G_{hom}(u_n,\tau_n) = G_{hom}(u,\tau).
\]
For each $n$, we have shown the existence of $(u_n^\e,\tau_n^\e)$ such that $u_n^\e \to u_n$ and  $\tau_n^\e \to \tau_n$ in $L^p(\Omega)$ and $L^q(\Omega; \bR^N)$ respectively and
\[
\limsup_{\e \to 0}  G_\e (u_n^\e,\tau_n^\e) \le G_{hom}(u_n,\tau_n),
\]
so we get
\[
\limsup_n \limsup_{\e \to 0}  G_\e (u_n^\e,\tau_n^\e) \le G_{hom}(u,\tau).
\]
From the same diagonalization argument  as before, there exists some $n(\e)$ such that $n(\e) \to \infty$ as $\e \to 0$  for which, by defining $(u_\e,\tau_\e) : = (u_{n(\e)}^\e,\tau_{n(\e)}^\e)$ we obtain
\[
\hbox{$u_\e \to u$ strongly in $L^p(\Omega)$, 
 $\tau_\e \to \tau$ strongly  in $L^q(\Omega; \bR^N)$}
 \]
and
\[
\limsup_{\e \to 0}  G_\e (u_n^\e,\tau_n^\e) \le G_{hom}(u,\tau).
\]
 This concludes  the proof of Lemma \ref{proof1}.\hfill $\square$\\

\textbf{Proof of Lemma  \ref{proof2}.}
Let $(u,\tau) \in W^{1,p}(\Omega) \times L^q(\Omega; \bR^N)$ and $f \in L^q(\Omega; \bR^N)$ with ${\rm div} (f)=0.$ We  assume  that $u_\e \to u$ strongly in $L^p(\Omega)$ and $\tau_\e \to \tau_\e$ in ${\mathcal T}.$
 For constant vectors $a_i,b_i,\eta_i \in \bR^n$, consider as before functions $\tilde{\phi_i} \in W^{1,p}_\#(Q)$ and $\tilde{g_i} \in L^q_\#(Q; \bR^N)$ such that
\[
L_{hom}(a_i,b_i + \eta_i) = \frac{1}{|Q|}\int_Q L\big( x,a_i + \nabla \tilde{\phi}_i(x),b_i + \eta_i + \tilde{g}_i(x) \big) \dx.
\]
 Denote $\partial_1 L$ the subdifferential of $L$ with respect to the middle variable and $\partial_2 L$ the subdifferential of $L$ with respect to the last variable. From the above we have both
\begin{equation}\label{min1}
{\rm div}\Big(\partial_1 L\big( y,a_i + D\tilde{\phi}_i(y),b_i + \eta_i + \tilde{g}_i(y) \big) \Big) = 0\, \,  a.e. \,\,  y \in Q,
\end{equation}
and
\begin{equation}\label{min2}
\int_Q \langle \partial_2 L\big( y,a_i + D\tilde{\phi}_i(y),b_i + \eta_i + \tilde{g}_i(y) \big),g(y)\rangle \dy = 0, 
\end{equation}
for any $g \in L^q_\#(Q; \bR^N)$. It follows from (\ref{min2}) that
\begin{equation}\label{min0}
 \partial_2 L\big( y,a_i + \nabla \tilde{\phi}_i(y),b_i + \eta_i + \tilde{g}_i(y) \big)=\nabla w(y)\qquad \text{ a.e. } y \in Q,
\end{equation}
for some $w \in W^{1,p}_\#(Q).$  It also follows from  Lemma \ref{ext1} that $w$ can be extended by periodicity to an element in $W_{loc}^{1,p}(\bR^N).$
 Now, let $\hat{u} \in W^{1,p} (\Omega) $ be a piecewise affine functions and $\hat{\tau}\in L^q(\Omega; \bR^N)$ be a piecewise constant function such that for some partition $\{ \Omega_i \}$ of $\Omega$ we have
\[
\hbox{$\hat{u}(x) = \langle a_i,x \rangle + \alpha_i \mbox{ for } x \in \Omega_i$ and  $
\hat{\tau}(x) =  b_i \mbox{ for } x \in \Omega_i.$}
\]
Consider now for $x \in \Omega_i$,
\[
\hbox{$\hat{u}_\e(x) := \hat{u}(x) + \e \tilde{\phi}_i(\frac{x}{\e})$ 
and $\hat{\tau}_\e(x) : = \hat{\tau}(x).$}
\]
From the convexity of $L$ we get
\begin{eqnarray*}
L\big(\frac{x}{\e},\nabla u_\e(x),\tau _\e(x) + f(x) \big) & \ge & L\big(\frac{x}{\e},\nabla \hat{u}_\e(x),\hat{\tau}_\e(x) + \eta_i + \tilde{g}_i(\frac{x}{\e}) \big)  \\
 & & +\langle \partial_1 L\big(\frac{x}{\e},\nabla \hat{u}_\e(x),\hat{\tau}_\e(x) + \eta_i + \tilde{g}_i(\frac{x}{\e}) \big) , \nabla u_\e(x) - \nabla \hat{u}_\e(x) \rangle \\
 & &+ \langle \partial_2 L\big(\frac{x}{\e},\nabla \hat{u}_\e(x),\hat{\tau}_\e(x) + \eta_i + \tilde{g}_i(\frac{x}{\e}) \big) , \tau_\e(x) - \hat{\tau}_\e(x) \rangle \\
 &  &+ \langle \partial_2 L\big(\frac{x}{\e},\nabla \hat{u}_\e(x),\hat{\tau}_\e(x) + \eta_i + \tilde{g}_i(\frac{x}{\e}) \big) , f(x) - \eta_i - \tilde{g}_i(\frac{x}{\e}) \rangle.
\end{eqnarray*}

  Consider now smooth functions $\Psi_i : \Omega_i \to \bR$ with compact support such that $0 < \Psi_i < 1$. Multiplying the above convexity inequality by $\Psi_i$, integrating over $\Omega_i$ and adding over all $i$, we get the following:
\begin{eqnarray*}
\S L(\frac{x}{\e},\nabla u_\e,\tau_\e + f)\dx &\ge&
\sum_i \int_{\Omega_i}  L\big(\frac{x}{\e},a_i + \nabla \tilde{\phi}_i(\frac{x}{\e}),b_i + \eta_i + \tilde{g}_i(\frac{x}{\e}) \big) \Psi_i(x)   \dx
\\&&
+ \sum_i \int_{\Omega_i} \langle \partial_1 L\big(\frac{x}{\e},a_i + \nabla \tilde{\phi}_i(\frac{x}{\e}),b_i + \eta_i + \tilde{g}_i(\frac{x}{\e}) \big) , \nabla u_\e(x) - \nabla \hat{u}_\e(x) \rangle\Psi_i(x) \dx
\\&&
+ \sum_i \int_{\Omega_i} \langle \partial_2 L\big(\frac{x}{\e},a_i + \nabla \tilde{\phi}_i(\frac{x}{\e}),b_i + \eta_i + \tilde{g}_i(\frac{x}{\e}) \big) , \tau_\e(x) - \hat{\tau}_\e(x) \rangle\Psi_i(x) \dx
\\&&
 + \sum_i \int_{\Omega_i} \langle \partial_2 L\big(\frac{x}{\e},a_i + \nabla \tilde{\phi}_i(\frac{x}{\e}),b_i + \eta_i + \tilde{g}_i(\frac{x}{\e}) \big) , f(x) - \eta_i \rangle\Psi_i(x) \dx \\&&
 + \sum_i \int_{\Omega_i} \langle \partial_2 L\big(\frac{x}{\e},a_i + \nabla \tilde{\phi}_i(\frac{x}{\e}),b_i + \eta_i + \tilde{g}_i(\frac{x}{\e}) \big) , - \tilde{g}_i(\frac{x}{\e}) \rangle\Psi_i(x) \dx .
\end{eqnarray*}

 Now we deal with each term independently. For the first term on the right hand side of the above expression we have
\[
 \int_{\Omega_i} L\big(\frac{x}{\e},a_i + \nabla \tilde{\phi}_i(\frac{x}{\e}),b_i + \eta_i + \tilde{g}_i(\frac{x}{\e}) \big)\Psi_i(x) \dx  \to  \int_{\Omega_i} L_{hom}(a_i,b_i + \eta_i)\Psi_i(x) \dx,
\]
by virtue of Lemma \ref{Riem-Leb}.

 For the second term, by integrating by parts and by then  taking into account \eqref{min1} we  obtain
\begin{eqnarray*}
\int_{\Omega_i} \langle \partial_1 L\big(\frac{x}{\e},a_i + \nabla \tilde{\phi}_i(\frac{x}{\e}),b_i + \eta_i + \tilde{g}_i(\frac{x}{\e}) \big) , \nabla u_\e(x) - \nabla \hat{u}_\e(x) \rangle\Psi_i(x) \dx \\
= -\int_{\Omega_i} \langle \partial_1 L\big(\frac{x}{\e},a_i + \nabla \tilde{\phi}_i(\frac{x}{\e}),b_i + \eta_i + \tilde{g}_i(\frac{x}{\e}) \big) , (u_\e - \hat{u}_\e)\nabla\Psi_i(x) \rangle \dx.
\end{eqnarray*}
It follows from Lemma \ref{Riem-Leb} and Proposition \ref{SubDifAve} below, that if $\e \to 0$ then,  
\begin{equation*}
\int_{\Omega_i} \langle \partial_1 L\big(\frac{x}{\e},a_i + \nabla \tilde{\phi}_i(\frac{x}{\e}),b_i + \eta_i + \tilde{g}_i(\frac{x}{\e}) \big) , (u_\e - \hat{u}_\e)\nabla\Psi_i(x) \rangle \dx \to\int_{\Omega_i}\langle \partial_1 L_{hom}(a_i,b_i+\eta_i), (u - \hat{u})\nabla\Psi_i(x) \rangle \dx
\end{equation*}
Integrate by parts on more time to get
\begin{equation*}
\int_{\Omega_i}\langle \partial_1 L_{hom}(a_i,b_i+\eta_i), (u - \hat{u})\nabla\Psi_i(x) \rangle \dx=-\int_{\Omega_i}\langle \partial_1 L_{hom}(a_i,b_i+\eta_i), \nabla u- \nabla \hat{u}\rangle \Psi_i(x) \dx,
\end{equation*}
from which one has
\begin{eqnarray*}
\int_{\Omega_i} \langle \partial_1 L\big(\frac{x}{\e},a_i + \nabla\tilde{\phi}_i(\frac{x}{\e}),b_i + \eta_i + \tilde{g}_i(\frac{x}{\e}) \big) , \nabla u_\e(x) - \nabla \hat{u}_\e(x) \rangle\Psi_i(x) \dx \to \\
\int_{\Omega_i}\langle \partial_1 L_{hom}(a_i,b_i+\eta_i), \nabla u - \nabla \hat{u} \rangle \Psi_i(x)\dx. 
\end{eqnarray*}
 For the third term, we use  (\ref{min0}) to get  $\partial_2 L\big(\frac{x}{\e},a_i + D\tilde{\phi}_i(\frac{x}{\e}),b_i + \eta_i + \tilde{g}_i(\frac{x}{\e}) \big)= \nabla w(\frac{x}{\e})$ for some $w \in W^{1,p}_{\#}(Q).$ Using an integration by parts,  we obtain
\begin{eqnarray*}
\int_{\Omega_i} \langle \nabla w(\frac{x}{\e}) , \tau_\e(x) - \hat{\tau}_\e(x) \rangle\Psi_i(x) \dx
&=&
-\int_{\Omega_i}\e w(\frac{x}{\e}) {\rm div} \big (\tau_\e(x)-\hat{\tau}_\e(x)\big ) \rangle\Psi_i(x) \dx\\&& -\int_{\Omega_i}\e w(\frac{x}{\e})  \langle \nabla \Psi_i(x) , \tau_\e(x) - \hat{\tau}_\e(x) \rangle \dx, 
\end{eqnarray*}
which goes to $0$ as $\e \to 0$ since $\tau_\e \to \tau $  in the ${\mathcal T}$-topology.

Similarly as above, the fourth term can be seen to converge to
\[
\int_{\Omega_i} \langle \partial_2 L_{hom}(a_i,b_i+\eta_i)  , f(x)- \eta_i\rangle\Psi_i(x) \dx, 
\]
while for the fifth term, we first observe that the function
\[
m_i(x):= \langle \partial_2 L\big(x,a_i + \nabla \tilde{\phi}_i(x),b_i + \eta_i + \tilde{g}_i(x) \big) , \tilde{g}_i(x)\rangle
\]
is $Q$-periodic, and thus setting $(m_i)_\e(x) := m_i(\frac{x}{\e})$, it follows from Lemma  \ref{Riem-Leb} that $
(m_i)_\e \rightharpoonup \overline{m}_i$ 
weakly in $L^1$, where
\[
\overline{m}_i =
\frac{1}{|Q|} \int_Q\langle \partial_2 L\big(y,a_i + \nabla \tilde{\phi}_i(y),b_i + \eta_i + \tilde{g}_i(y) \big) , - \tilde{g}_i(y) \rangle \dy, 
\]
which is equal to $0$ in view of \eqref{min2}.
The fifth term therefore  disappears as $\e \to 0$. 

Putting now 
all of the above together we obtain that 
\begin{eqnarray*}
\liminf_{\e \to 0} \S L(\frac{x}{\e},\nabla u_\e,\tau_\e + f ) \dx &\ge&
 \sum_i \int_{\Omega_i} L_{hom}( a_i,b_i+\eta_i ) \Psi_i(x)\dx\\
&& + \sum_i \int_{\Omega_i}\langle \partial_1 L_{hom}( a_i,b_i+\eta_i )  , \nabla u(x) - \nabla \hat{u}(x) \rangle\Psi_i(x)\dx \\
&& + \sum_i \int_{\Omega_i}\langle \partial_2 L_{hom}( a_i,b_i+\eta_i )  , f(x) - \eta_i \rangle\Psi_i(x) \dx.
\end{eqnarray*}
By taking into account the estimate
\[|\partial L_{hom} (a,b)| \le M (1+|a|^{p-1}+|b|^{q-1})  \qquad \text{ for all } a,b \in \bR^N,\]
which follows from estimate (\ref{bd-hom}) in  Proposition \ref{l-infty},  and letting  $\Psi_i\uparrow 1$ on each  $\Omega_i,$ it follows from the dominated convergence theorem that
\begin{eqnarray*}
\liminf_{\e \to 0} \S L(\frac{x}{\e},\nabla u_\e,\tau_\e + f ) \dx &\ge&
 \S L_{hom}\big( \nabla \hat{u}(x),\hat{\tau}(x) + \tilde{f}(x) \big) \dx\\&&
 + \S \langle \partial_1 L_{hom}\big( \nabla \hat{u}(x),\hat{\tau}(x) + \tilde{f}(x) \big) , \nabla u(x) - \nabla \hat{u}(x) \rangle \dx \\
&& + \S \langle \partial_2 L_{hom}\big( \nabla \hat{u}(x),\hat{\tau}(x) + \tilde{f}(x) \big) , f(x) - \tilde{f}(x) \rangle \dx.
\end{eqnarray*}
where $\tilde{f} \in L^q(\Omega; \bR^N)$ is a function defined by $\tilde{f}(x)=\eta_i$ on $\Omega_i.$
 The above is valid for arbitrary  piecewise affine function $\hat{u}$, and  piecewise constant functions $\hat{\tau}, \tilde{f}$. We can then let $\hat{u} \to u$ in $W^{1,p}(\Omega)$ and $\hat{\tau} \to \tau$  and $\tilde{f} \to f$ in $L^q(\Omega; \bR^N)$ to obtain
\[
\liminf_{\e \to 0} \S L(\frac{x}{\e},\nabla u_\e,\tau_\e + f ) \dx \ge \S L_{hom}\big( \nabla u(x),\tau(x) + f(x) \big) \dx.
\]
This completes the proof.
\hfill $\square$

Before proceeding to the next subsection, we note the following slight  extension of Lemma  \ref{proof1}, which will be needed for Proposition \ref{gam-self} below. We note  that the proof    is  known when $G_\e$ is independent  of the second variable, and here we show that the same proof with minor modification works for general Lagrangians just as in Theorem \ref{gam-lim}.
\begin{lem}\label{boundary}
Let $G_\e$ and $G_{hom}$ be as in Theorem \ref{gam-lim}. Then, for any $(u,\tau) \in W^{1,p}(\Omega) \times L^q(\Omega;\bR^N)$, there exist a sequence $(u_\e,\tau_\e)$ such that $u - u_\e  \rightharpoonup 0$ weakly in $W^{1,p}(\Omega)$ and $\tau_\e \to \tau$ in the ${\mathcal T}$-topology. Furthermore, $u - u_\e \in W^{1,p}_0(\Omega)$ and for this sequence:
\[
\limsup_{\e \to 0}  G_\e (u_\e , \tau_\e ) \le G (u,\tau).
\]
\end{lem}
\textbf{Proof.}
From Theorem \ref{gam-lim}, there exist a sequence $(\tilde{u}_\e,\tau_\e)$ with $\tilde{u}_\e \to u$ in $L^p( \Omega )$ and $\tau_\e \to \tau$ in the ${\mathcal T}$-topology, such that
\[
G_{hom}(u,\tau) = \lim_{\e \to 0}  G_\e(\tilde{u}_\e,\tau_\e).
\]
 Up to a subsequence one may assume that
\[
\tilde{u}_\e \rightharpoonup u \mbox{ weakly in } W^{1,p}(\Omega).
\]
 Pick any $\phi \in W^{1,p}_0( \Omega )$ with $\phi > 0$ in $\Omega$. Define:
\[
u_\e(x) : = \left\{
\begin{array}{lr}
\tilde{u}_\e & u(x) - \phi(x) \le \tilde{u}_\e \le u(x) + \phi(x) \\
u(x) - \phi(x) & \tilde{u}_\e (x) < u(x) - \phi(x) \\
u(x) + \phi(x) & u(x) + \phi(x) < \tilde{u}_\e(x)
\end{array}
\right. .
\]
 Note that $u_\e-u \in W^{1,p}_0( \Omega )$ and since $\tilde{u}_\e \rightharpoonup u$ weakly in $W^{1,p}(\Omega)$, so $u_\e \rightharpoonup u$ weakly in $W^{1,p}(\Omega)$.
Note that  $L+ C_0 \ge 0.$ For any $f \in L^q(\Omega;\bR^N)$ with ${\rm div}(f) = 0$ we have
\begin{eqnarray*}
G_\e (u_\e,\tau_\e) + C_0| \Omega | &\le&
 \int_{ \{ u_\e \ne \tilde{u}_\e \} } \Big[ L \big( \frac{x}{\e},\nabla u_\e(x),f(x) + \tau_\e(x) \big) + C_0 \Big]\, dx \\&&+
 \int_{ \{u_\e = \tilde{u}_\e \} } \Big[ L \big( \frac{x}{\e},\nabla u_\e(x),f(x) + \tau_\e(x) \big) + C_0 \Big]\, dx.
\end{eqnarray*}

 For $x$ in the set $\{ u_\e \ne \tilde{u}_\e \}$, the norm of $\nabla \tilde{u}_\e(x)$ is  controlled by the norm of  $|\nabla u(x)| + |\nabla \phi (x)|$. It follows that
\begin{eqnarray*}
G_\e (u_\e,\tau_\e) + C_0| \Omega | &\le
& \int_{ \{ u_\e \ne \tilde{u}_\e \} } \Big[ C_1\Big ((|\nabla u(x)| + |\nabla \phi (x)|)^p + |\tau_\e(x)+f(x)|^q + 1 \Big) + C_0 \Big]\, dx \\&&+
 \S \Big[ L \big( \frac{x}{\e},\nabla \tilde{u}_\e(x),f(x) + \tau_\e(x) \big) + C_0 \Big]\,dx.
\end{eqnarray*}
 Take now the infimum over all  $f \in L^q(\Omega;\bR^N)$ with ${\rm div}(f) = 0$ and subtract the latter by $C_0|\Omega |$. Since $|\{ u_\e \ne \tilde{u}_\e \}| \to 0$ and $G_{hom}(u,\tau) = \lim_\e G_\e(\tilde{u}_\e,\tau_\e)$, we obtain
\[
\limsup_{\e \to 0}  G_\e (u_\e,\tau_\e) \le G(u,\tau).
\]

\subsection{Variational homogenization  of   maximal monotone operators on $W_0^{1,p}(\Omega)$} \label{CMMO}
In this section we establish a homogenization result for selfdual Lagrangians  on   $W_0^{1,p}(\Omega) \times W^{-1,q}(\Omega)$ where $\frac{1}{p}+\frac{1}{q}=1$ and then proceed to prove Theorem \ref{main}.

\begin{thm}\label{gam-mm} Let $\Omega$  be a regular bounded domain, 
$Q$ be an open non-degenerate parallelogram in  $\bR^n$, and $L: \Omega \times \mathbb{R}^{N} \times \mathbb{R}^{N}\to \R $ be  an  $\Omega$-dependent selfdual Lagrangian such that:\\
(1) For each $a,b \in \bR^N $ the function
$x \to L(x,a,b)$  is $Q$-periodic,\\
(2)  For some constants $C_0,C_1 \ge 0$, we have for every $x\in \bR^N$,
\begin{equation} \label{es300}
C_0(|a|^p + |b|^q ) \le L(x,a,b) \le C_1(|a|^p + |b|^q + 1), 
\end{equation}
where $p>1$ and $\frac{1}{p} + \frac{1}{q} = 1$. Let $u^*_n \to u^* $ strongly in $W^{-1,q}(\Omega)$ and let   $u_n$ be solutions and $\tau_n$ be momenta for the Dirichlet boundary value problems
\begin{eqnarray}\label{mm}
\left\{
  \begin{array}{lll}
    \tau_n(x) \in \bar \partial L(\frac{x}{\e_n}, \nabla u_n(x)) & a.e. \quad x \in \Omega \\
   -{\rm div} (\tau_n (x))= u^*_n(x) & x \in \Omega\\
u_n \in W_0^{1,p}(\Omega).\\
  \end{array}
\right.
\end{eqnarray}
Then, up to a subsequence, 
\[
\hbox{$u_n \to u$ weakly in  $W_0^{1,p}(\Omega)$ and  $\tau_n \to \tau$ weakly in $L^q(\Omega; \bR^N)$,}
\]
where $u$ is a solution and $\tau$ is a momentum of the homogenized problem
\begin{eqnarray}\label{mm-hom}
\left\{
  \begin{array}{lll}
    \tau(x) \in \bar \partial L_{hom}( \nabla u(x)) & a.e. \quad  x \in \Omega \\
   -{\rm div} (\tau (x))= u^*(x) & x \in \Omega\\
u \in W_0^{1,p}(\Omega), \\
  \end{array}
\right.
\end{eqnarray}
where $L_{hom}$ is the selfdual Lagrangian on $\bR^N \times \bR^N$
defined by
\begin{equation}\label{hom301}
L_{hom}(a,b) := \min_{\substack{\phi \in W^{1,p}_\# (Q) \\ g \in L_\#^q(Q; \bR^N)}}\frac{1}{|Q|}\int_Q L\big( x,a + D\phi(x),b + g(x) \big) \dx.
\end{equation}
\end{thm}
This will follow from the following proposition. 

\begin{pro}\label{gam-self} Let $\Omega, Q$ and 
$L$
be  
as in Theorem \ref{gam-mm}, and let $\{F_\e; \e >0\}$ be the family of selfdual Lagrangians on  $W_0^{1,p}(\Omega) \times W^{-1,q}(\Omega)$ defined by
\[
F_\e (u,u^*) := \inf_{\substack{f \in L^q(\Omega; \bR^N) \\ -{\rm div} f = u^*}} \S L\big(\frac{x}{\e},\nabla u(x), f(x)\big) \dx.
\]
Then, there exists a selfdual Lagrangian  $F_{hom}$ on $\bR^N \times \bR^N$ that is a $\Gamma$-limit  of $\{F_\e; \e >0\}$ on $W_0^{1,p} (\Omega) \times W^{-1,q}(\Omega)$.  It is given by the formula 
\[
F_{hom} (u,u^*) := \inf_{\substack{f \in L^q(\Omega; \bR^N) \\ -{\rm div} f = u^*}} \S L_{hom}\big(\nabla u(x), f(x)\big) \dx,
\]
where $L_{hom}$ is the selfdual Lagrangian on $\bR^N \times \bR^N$ defined by (\ref{hom301}),  
and which satisfies for all $(a,b)\in \bR^N \times \bR^N$ 
\[
 C_0(|a|^p + |b|^q -1) \le L_{hom}(a,b) \le C_1(|a|^p + |b|^q + 1).
\]
\end{pro}
\textbf{Proof.} Note first that the selfduality and uniform bounds of $L_{hom}$ follow from Proposition \ref{l-infty}. It also follows from Proposition \ref{lift} that both $F_\e$ and $F_{hom}$ are selfdual Lagrangians on $W_0^{1,p}(\Omega) \times W^{-1,q}(\Omega).$  Given $(u,u^*) \in W_0^{1,p}(\Omega) \times W^{-1,q}(\Omega)$, we now show the existence of a sequence $\{(u_\e, u^*_\e) \in W_0^{1,p}(\Omega) \times W^{-1,q}(\Omega)\}$ with  $u_\e \to u$ weakly in $W_0^{1,p}(\Omega)$ and $u^*_\e \to u^*$ strongly in $W^{-1,q}(\Omega)$ and such that 
\begin{equation}\label{m-sup}
\limsup_{\e \to 0}  F_{\e} (u_{\e},u^*_{\e}) \le F_{hom}(u,u^*).
\end{equation}
For that we consider  $\{G_\e; \e >0\}$ be a family of functionals on  $W^{1,p}(\Omega) \times L^q(\Omega;\bR^N)$ defined by
\[
G_\e (u,\tau) := \inf_{\substack{f \in L^q(\Omega;\bR^N) \\ {\rm div} f = 0}} \S L\big(\frac{x}{\e},\nabla u(x),\tau(x) + f(x)\big) \dx, 
\]
and
\[
G_{hom} (u,\tau) := \inf_{\substack{f \in L^q(\Omega;\bR^N) \\ {\rm div} f = 0}} \S L_{hom}\big(\nabla u(x),\tau(x) + f(x)\big) \dx,
\]
 Take $\tau \in L^q(\Omega; \bR^N)$ such that ${\rm div}(\tau)=u^*.$ It follows from Lemma  \ref{proof1} and Lemma 
\ref{boundary} that there exists $(u_\e, \tau_\e) \in W_0^{1,p}(\Omega) \times L^q(\Omega; \bR^N)$ such that  $u_\e \to u$ strongly in $L^p(\Omega)$ and $\tau_\e \to \tau $ strongly in $L^q(\Omega; \bR^N)$ and
\[
\limsup_{\e \to 0} G_\e (u_{\e},\tau_{\e}) \le  G_{hom}(u,\tau).
\]
The sequence $u_\e$ is bounded in $W_0^{1,p}(\Omega)$, so we may assume $u_\e \to u$ weakly in $W_0^{1,p}(\Omega).$
Since $\tau_\e\to \tau $ strongly in $L^q(\Omega; \bR^N),$ it follows that $u^*_\e :={\rm div}(\tau_\e)\to {\rm div}(\tau)=u^*$ strongly in $W^{-1,q}(\Omega).$ Thus, the inequality (\ref{m-sup}) follows by noticing that
$G_\e  (u_{\e},\tau_{\e})=F_{\e} (u_{\e},u^*_\e)$ and $G_{hom} (u,\tau)=F_{hom} (u,u^*).$\\

We shall now show that if  $(u,u^*) \in W_0^{1,p}(\Omega) \times W^{-1,q}(\Omega)$ and   $u_\e \to u$ weakly in $W_0^{1,p}(\Omega)$ and $u^*_\e\to u^*$ strongly in $W^{-1,q}(\Omega)$ then
\begin{equation}\label{m-inf}
F_{hom}(u,u^*) \le \liminf_{\e \to 0}  F_{\e} (u_{\e},u^*_{\e}).
\end{equation}
 Take an arbitrary element in $(v,v^*) \in W_0^{1,p}(\Omega) \times W^{-1,q}(\Omega).$ From the above,  there exists   $(v_\e, v^*_\e) \in W_0^{1,p}(\Omega) \times W^{-1,q}(\Omega)$ with  $\tilde v_\e \to v$ weakly in $W_0^{1,p}(\Omega)$ and $v^*_\e \to v^*$ strongly in $W^{-1,q}(\Omega)$ and such that 
\begin{equation*}\label{m-sup.0}
\limsup_{\e \to 0}  F_{\e} (v_{\e},v^*_{\e}) \le F_{hom}(v,v^*).
\end{equation*}
By the self duality of $F_\e $ we have
\begin{eqnarray*}
F_{\e} (u_{\e},u^*_{\e})=F^*_{\e} (u^*_{\e}, u_{\e}) &= & \sup\{ \langle u_{\e}, w^*\rangle +\langle u^*_{\e}, w\rangle - F_{\e}(w,w^*); (w,w^*) \in W_0^{1,p}(\Omega) \times W^{-1,q}(\Omega)  \}\\ &\ge & \langle u_{\e}, v^*_\e\rangle +\langle u^*_{\e}, v_\e\rangle - F_{\e}(v_\e,v^*_\e), 
\end{eqnarray*}
from which we get
\begin{eqnarray*}
\liminf_{\e \to 0} F_{\e} (u_{\e},u^*_{\e})  &\ge & \liminf_{\e \to 0}  \{\langle u_{\e}, v^*_\e\rangle +\langle u^*_{\e}, v_\e\rangle - F_{\e}(v_\e,v^*_\e) \}\\
&= & \langle u, v^*\rangle +\langle u^*, v \rangle-\limsup_{\e \to 0}  F_{\e}(v_\e,v^*_\e) \\
&\ge & \langle u, v^*\rangle +\langle u^*, v \rangle-F_{hom}(v,v^*).
\end{eqnarray*}
Since the above holds for an arbitrary $(v,v^*) \in W_0^{1,p}(\Omega) \times W^{-1,q}(\Omega),$ we obtain
\begin{eqnarray*}
F^*_{hom}(u^*,u) \le \liminf_{\e \to 0} F_{\e} (u_{\e},u^*_{\e}).
\end{eqnarray*}
Taking into consideration that $F_{hom}$ is selfdual we obtain
\begin{eqnarray*}
F_{hom}(u,u^*)=F^*_{hom}(u^*,u) \le \liminf_{\e \to 0} F_{\e} (u_{\e},u^*_{\e}), 
\end{eqnarray*}
as desired.
 \hfill $\square$\\

\textbf{Proof of Theorem \ref{gam-mm}.} Since    
$(u_n, \tau_n)$ are solutions of  (\ref{mm}), it follows that 
\begin{eqnarray}\label{mm1}
0&=&\int_{\Omega} L(\frac{x}{\e_n},\nabla u_n(x), \tau_n (x) ) \dx-\int_{\Omega} \langle \nabla u_n(x),\tau_n(x)\rangle_{\bR^N}\,dx \nonumber \\&=&\int_{\Omega} L(\frac{x}{\e_n},\nabla u_n(x), \tau_n (x)) \, dx-\int_{\Omega}  u_n(x)u^*_n(x)\,dx. 
\end{eqnarray}
Due to the coercivity assumption on $L$ and the strong convergence of $u^*_n,$ it follows that the sequence $u_n$ is bounded in $W_0^{1,p}(\Omega)$ and $\tau_n$ is bounded in $L^q(\Omega; \bR^N).$ Thus, up to a subsequence,  $u_n \to u $  weakly in  $ W_0^{1,p}(\Omega)$ and $\tau_n \to \tau $  weakly in $L^q(\Omega; \bR^N).$ We also have ${\rm div}(\tau_n)=u^*_n \to u^*= {\rm div} (\tau)$ strongly in $W^{-1,q}(\Omega),$ from which we indeed have $\tau_n \to \tau $   in the ${\mathcal T}$-topology (introduced in Theorem \ref{gam-lim}).

Taking $f \in L^q(\Omega; \bR^N)$ with $ {\rm div} f = 0 $, it follows from (\ref{mm1}) that
 \begin{eqnarray*}
\int_{\Omega} L(\frac{x}{\e_n},\nabla u_n(x), \tau_n (x) ) \dx&=&\int_{\Omega} u_n(x)u^*_n(x)\,dx\\
&=&-\int_{\Omega}  u_n(x){\rm div}(\tau_n+f)\,dx\\
&=&\int_{\Omega}  \langle \nabla u_n(x), \tau_n+ f\rangle_{\bR^N} \,dx\\
&\le & \S L\big(\frac{x}{\e_n},\nabla u_n(x), \tau_n +f(x)\big) \, dx.
\end{eqnarray*}

This indeed shows that
\[
\S L\big(\frac{x}{\e_n},\nabla u_n(x), \tau_n(x)\big) \dx=\inf_{\substack{f \in L^q(\Omega; \bR^N)\\ {\rm div} f = 0}} \S L\big(\frac{x}{\e_n},\nabla u_n(x), \tau_n(x)+f(x)\big) \dx.
\]
Let \[
G_{\e_n} (v,\hat \tau) := \inf_{\substack{f \in L^q(\Omega;\bR^N) \\ {\rm div} f = 0}} \S L\big(\frac{x}{\e_n},\nabla v(x),\hat \tau(x) + f(x)\big) \dx.
\]
It then follows that $\S L\big(\frac{x}{\e_n},\nabla u_n(x), \tau_n(x)\big) \dx=G_{\e_n}(u_n, \tau_n).$
Define $H:W_0^{1,p}(\Omega) \times L^q(\Omega; \bR^N)\to \bR $ by $H(v,\tilde \tau )= \int_{\Omega} \langle \nabla v(x), \tilde \tau(x)\rangle_{\bR^N}\, dx.$ Note that $H$ is continuous if we consider the weak topology of  $W_0^{1,p}(\Omega)$ and the ${\mathcal T} $-topopogy for  $L^q(\Omega; \bR^N)$. It then follows from Lemma  \ref{proof2} that
\begin{eqnarray*}
\S L_{hom}\big(\nabla u(x), \tau(x)\big) \dx-H(u,\tau)  &\le & \liminf_{\e_n\to 0} \big [G_{\e_n}(u_n,\tau_n)-H(u_n,\tau_n) \big ]\\
 &=&\liminf_{\e_n\to 0}\Big [\int_{\Omega} L(\frac{x}{\e_n},\nabla u_n(x), \tau_n (x) ) \dx-\int_{\Omega} u_n(x){\rm div}(\tau_n(x))\,dx\Big ]
\\&=& 0.
\end{eqnarray*}
On the other hand, we have that 
\begin{eqnarray}
\S L_{hom}\big(\nabla u(x), \tau(x)\big) \dx-H(u,\tau)  = \S \Big [L_{hom}\big(\nabla u(x), \tau(x)\big)-\langle \nabla  u(x),\tau(x)\rangle_{\bR^N}\Big ] \,dx\ge0.
\end{eqnarray}
which means that the latter is indeed zero, i.e.,
\[\S \Big [L_{hom}\big(\nabla u(x), \tau(x)\big)-\langle \nabla  u(x),\tau(x)\rangle_{\bR^N} \,\Big] dx=0.\]
Since the integrand is itself non-negative we have
\[L_{hom}\big(\nabla u(x), \tau(x)\big)-\langle \nabla  u(x),\tau(x)\rangle_{\bR^N}=0 \quad \text{ a.e. } x \in \Omega,\]
which together with $-{\rm div} (\tau (x))= u^*(x)$, yields
\begin{eqnarray*}
\left\{
  \begin{array}{lll}
    \tau(x) \in \bar \partial L_{hom}( \nabla u(x)), & a.e. x \in \Omega, \\
   -{\rm div} (\tau (x))= u^*(x), & x \in \Omega,\\
u \in W_0^{1,p}(\Omega).\\
  \end{array}
\right.
\end{eqnarray*}

\section{Appendix}

We shall here state some of the results  used throughout the proof.

\begin{lem}\label{II}
Assume $L:\bR^n\times\bR^n\to \bR$ is a convex function such that $C_0(|a|^p + |b|^q - 1) \le L(a,b) \le C_1(|a|^p + |b|^q + 1)$ for all $a,b \in \bR^N$ where $p,q>1$ are two constants. Suppose $\Omega$ is a bounded open domain in $\bR^N$ and   $\tau_1 \in L^p(\Omega; \bR^N)$ and $\tau_2 \in L^q(\Omega; \bR^N)$ are two   piecewise constant functions such that
\[
\hbox{$\tau_1(x) =  a_i, \quad x \in {\Omega}_i,$
and $
 \tau_2(x) =  b_i, \quad x \in {\Omega}_i,$}
\]
where  $\{{\Omega}_i\}_{i \in I}$ is a  finite polyhedral partitions of $\Omega$, and $\{a_i\}_{i \in I}, \{b_i\}_{i \in I}$ are two sequences $\in \bR^N.$
Then
\begin{equation*}
\min_{ \substack{f \in L^q(\Omega;\bR^N) \\ {\rm div} f = 0} } \S L\big( \tau_1,\tau_2(x)+f(x) \big) \dx \ge \sum_{i \in I} |\Omega_i| \inf_{\eta_i \in \bR^n} L( a_i,b_i + \eta_i).
\end{equation*}
\end{lem}

\textbf{Proof}
We first prove a stronger result (actually an equality) when  the set index $I$ is a singleton.   For any constant $ \eta\in \bR^N$ we have
\begin{eqnarray*}
\min_{ \substack{f \in  L^q(\Omega;\bR^N)  \\ {\rm div} f = 0} } \S L\big( a,b+f(x) \big) \dx  \le  \S L( a,b + \eta ) \dx  = |\Omega|L( a,b + \eta),
\end{eqnarray*}
from which we obtain
\begin{eqnarray*}
\min_{ \substack{f \in  L^q(\Omega;\bR^N)  \\ {\rm div} f = 0} } \S L\big( a,b+f(x) \big) \dx & \le & \inf_{\eta \in \bR^N}|\Omega|L( a,b + \eta),
\end{eqnarray*}
 Let now $\tilde{f}$ be the element in $ L^q(\Omega;\bR^N) $ with ${\rm div} \tilde{f} = 0$ such that
\[
\S L\big(a,b + \tilde{f}(x)\big) \dx = \min_{ \substack{f \in  L^q(\Omega;\bR^N) ) \\ {\rm div} f = 0} } \S L\big( a,b+f(x) \big) \dx.
\]
Using Jensen's inequality,  we obtain
\begin{eqnarray*}
\inf_{\eta \in \bR^N}|\Omega|L( a,b + \eta) & \le & |\Omega| L\big( a,b + \frac{1}{|\Omega|} \S \tilde{f}(x) \dx \big)\\& = & |\Omega| L\big( \frac{1}{|\Omega|} \S a \dx ,\frac{1}{|\Omega|} \S b + \tilde{f}(x) \dx \big)\\
 & \le & \S L\big(a,b + \tilde{f}(x)\big) \dx \\
 & = & \min_{ \substack{f \in  L^q(\Omega;\bR^N) \\ {\rm div} f = 0} } \S L\big( a,b+f(x) \big) \dx.
\end{eqnarray*}
This completes the proof for $I$ being a singleton. Now we prove it for the general case. Note first that, using  the above argument
 on each $\Omega_i$ we have
\begin{eqnarray}\label{div-equ1}
\inf_{\substack{g \in  L^q(\Omega_i;\bR^N)  \\ {\rm div} g = 0}} \int_{\Omega_i} L\big(a_i,b_i + g(x)\big) \dx = \inf_{\eta_i \in \bR^N}|\Omega_i|L( a_i,b_i + \eta_i).
\end{eqnarray}
One also can easily deduce that
\begin{eqnarray}\label{div-equ}
\inf_{\substack{f \in  L^q(\Omega;\bR^N)  \\ {\rm div} f = 0}} \S L\big(\tau_1(x),\tau_2(x) + f(x)\big) \dx \ge  \sum_i \inf_{\substack{f_i \in  L^q(\Omega_i;\bR^N)  \\ {\rm div} f_i = 0}} \int_{\Omega_i} L\big(a_i,b_i + f_i(x)\big) \dx.
\end{eqnarray}
In fact if $\inf_{\substack{f \in L^q(\Omega;\bR^N)  \\ {\rm div} f = 0}} \S L\big(\tau_1(x),\tau_2(x) + f(x)\big) \dx =\S L\big(\tau_1(x),\tau_2(x) + \bar f(x)\big) \dx$ for some $\bar f \in  L^q(\Omega;\bR^N) $ with ${\rm div} (\bar f)=0,$ then
\begin{eqnarray*}
\S L\big(\tau_1(x),\tau_2(x) + \bar f(x)\big) \dx&=&\sum_{i \in I}\int_{\Omega_i} L\big(a_i,b_i + \bar f(x)\big) \dx\\
& \geq& \sum_{i \in I}\inf_{\substack{f_i \in  L^q(\Omega_i;\bR^N) \\ {\rm div} f_i = 0}}\int_{\Omega_i} L\big(a_i,b_i + f_i(x)\big) \dx. 
\end{eqnarray*}
The proof therefore follows from combining (\ref{div-equ1}) and (\ref{div-equ}).
\hfill $\square$\\

 The following three Lemmas are standard and we refer to $\cite{suq}$ for the proof.

\begin{lem}\label{Riem-Leb}
Let $r \geq 1$ and  $f \in L^r(Q).$ Then $f$ can be extended by periodicity to a function (still denoted by $f$) belonging to $L_{loc}^r(\bR^N).$ Moreover, if $(\e_k)$ is a sequence of positive real numbers converging to $0$ and $g_k(x)=g(\frac{x}{\e_k})$.
\[\hbox{If $1 \leq r < \infty$, then $ f_k \to M(f)= \frac{1}{|Q|}\int_Q f (x) \dx$ weakly in  $L_{loc}^r(\bR^N)$,}
\]
and
\[
\hbox{if $r=\infty$, then   $f_k \to M(f)$ weak$^*$ in $L^{\infty}(\bR^N)$.}
\]
\end{lem}

\begin{lem}\label{ext1} Let $r>1$ and $u \in W_{\#}^{1,r}(Q)$, then $u$ can be extended by periodicity to an element of $W_{loc}^{1,r}(\bR^N).$
\end{lem}
\begin{lem}\label{ext2} Let $r>1$ and $r'=\frac{r}{r-1}.$ Let $ g \in L^{r'}(Q; \bR^N)$ such that $\int_{Q}\langle g(x), \nabla v(x) \rangle \, dx=0$ for every   $v \in W_{\#}^{1,r}(Q).$ Then $g$ can be extended by periodicity to an element of $L_{loc}^{r'}(\bR^N; \bR^N),$ still denoted by $g$,  such that ${\rm div} (g)=0$ in $D'(\bR^N).$
\end{lem}

\end{document}